\documentclass{amsproc}

\usepackage{amsmath}
\usepackage{amsfonts}
\usepackage{amssymb}
\usepackage{amscd}
\usepackage{verbatim}
 \usepackage[all]{xy}

\newtheorem{theorem}{Theorem}[section]
\newtheorem{lemma}[theorem]{Lemma}

\theoremstyle{definition}

\newtheorem{example}[theorem]{Example}
\newtheorem{proposition}[theorem]{Proposition}
\newtheorem{corollary}[theorem]{Corollary}

\theoremstyle{remark}
\newtheorem{remark}[theorem]{Remark}

\numberwithin{equation}{section}

\newcommand{\RR}{\mathbf{R}}
\newcommand{\NN}{\mathbf{N}}
\newcommand{\AAA}{\mathbf{A}}
\newcommand{\ZZ}{\mathbf{Z}}

\newcommand{\QQ}{\mathbf{Q}}

\newcommand{\OO}{\mathcal  {O}}

\newcommand{\fra}{\frak{a}}
\newcommand{\frb}{\frak{b}}

\newcommand{\frmm}{\frak{m}}

\DeclareMathOperator{\codim}{codim} 
\DeclareMathOperator{\Hom}{Hom}

 \DeclareMathOperator{\Spec}{Spec}

\DeclareMathOperator{\ord}{ord} \DeclareMathOperator{\mld}{mld}

\newcommand{\llbracket}{[\negthinspace[}
\newcommand{\rrbracket}{]\negthinspace]}

\newcommand{\llparenthesis}{(\negthinspace(}
\newcommand{\rrparenthesis}{)\negthinspace)}

\begin{document}

\title{Jet Schemes and Singularities}

\author{Lawrence Ein}
\address{Department of Mathematics \\ University of
Illinois at Chicago, 851 South Morgan Street (M/C 249)\\
Chicago, IL 60607-7045, USA} \email{ein@math.uic.edu}
\thanks{The first author was supported in part by NSF under Grant DMS-0200278.}

\author{Mircea Musta\c{t}\v{a}}
\address{Department of Mathematics
\\ University of Michigan \\ Ann Arbor, MI 48109, USA} \email{mustata@umich.edu}
\thanks{The second author was supported in part by NSF under Grants DMS-0500127
and DMS-0111298, and by a Packard Fellowship}

\maketitle

\section{Introduction}

The study of singularities of pairs is fundamental for higher
dimensional birational geometry. The usual approach to invariants of
such singularities is via divisorial valuations, as in
\cite{kollar}. In this paper we give a self-contained presentation
of an alternative approach, via contact loci in spaces of arcs. Our
main application is a version of Inversion of Adjunction for a
normal $\QQ$--Gorenstein variety embedded in a nonsingular variety.

The invariants we study are the minimal log discrepancies. Their
systematic study is due to Shokurov and Ambro, who made in
particular several conjectures, whose solution would imply the
remaining step in the Minimal Model Program, the Termination of
Flips (see \cite{ambro} and \cite{shokurov}). We work in the
following setting: we have a pair $(X,Y)$, where $X$ is a normal,
$\QQ$--Gorenstein variety and $Y$ is equal to a formal linear
combination $\sum_{i=1}^s q_iY_i$, where all $q_i$ are non-negative
real numbers, and all $Y_i$ are proper closed subschemes of $X$. To
every closed subset $W$ of $X$ one associates an invariant, the
minimal log discrepancy ${\rm mld}(W;X,Y)$, obtained by taking the
minimum of the so-called log discrepancies  of the pair $(X,Y)$ with
respect to all divisors $E$ over $X$ whose image lies in $W$. We do
not give here the precise definition, but refer instead to \S 7.

The space of arcs $J_{\infty}(X)$ of $X$ parametrizes morphisms
${\rm Spec}\,k\llbracket t\rrbracket\to X$, where $k$ is the ground field. It
consists of the $k$--valued points of a scheme that is in general not
of finite type over $k$. This space is studied  by looking at its
image  in the jet schemes of $X$ via the truncation maps. The
$m^{\rm th}$ jet scheme $J_m(X)$ is a scheme of finite type that
parametrizes morphisms ${\rm Spec}\,k[t]/(t^{m+1})\to X$. It was
shown in \cite{EMY} that the minimal log discrepancies can be
computed in terms of the codimensions of certain contact loci in
$J_{\infty}(X)$, defined by the order of vanishing along various
subschemes of $X$. As an application it was shown in \cite{EMY} and
\cite{EM} that a precise form of Inversion of Adjunction holds for
locally complete intersection varieties. In practice one always
works at the finite level, in a suitable jet scheme, and therefore
in order to apply the above-mentioned criterion one has to find (a
small number of) equations for the jets that can be lifted to the
space of arcs. This was the technical core of the argument in
\cite{EM}. In the present paper we simplify this approach by giving
first an interpretation of minimal log discrepancies in terms of the
dimensions of certain contact loci in the jet schemes, as opposed to
such loci in the space of arcs (see Theorem~\ref{new_interpretation}
for the precise statement). We apply this point of view to give a
proof of the following version of Inversion of Adjunction. This has
been proved independently also by Kawakita in \cite{kawakita2}.

\begin{theorem}\label{theorem_introduction}
Let $A$ be a nonsingular variety and $Y=\sum_{i=1}^sq_iY_i$, where
the $q_i$ are non-negative real numbers and the $Y_i$ are proper
closed subschemes of $A$. If $X$ is a closed normal subvariety of
$A$ of codimension $c$ such that $X$ is not contained in the support
of any $Y_i$, and if $rK_X$ is Cartier, then there is an ideal $J_r$
on $X$ whose support is the non-locally complete intersection locus
of $X$ such that
\begin{equation}\label{formula_introduction}
{\rm mld}(W;A,Y+cX)={\rm
mld}\left(W;X,Y\vert_X+\frac{1}{r}V(J_r)\right)
\end{equation}
 for every proper
closed subset $W$ of $X$.
\end{theorem}

\bigskip

When $X$ is locally complete intersection, this recovers the main
result from \cite{EM}. We want to emphasize that from the point of
view of jet schemes the ideal $J_r$ in the above theorem appears
quite naturally. In fact, the reduction to complete intersection
varieties is a constant feature in the the study of jet schemes
(see, for example, the results in \S 4). On the other hand, the
appearance of $\frac{1}{r}V(J_r)$ on the right-hand side of
(\ref{formula_introduction}) is the reason why the jet-theoretic
approach has failed so far to prove the general case of Inversion of
Adjunction.

The main ingredients in the arc-interpretation of the invariants of
singularities are the results of Denef and Loeser from \cite{DL}. In
particular, we use their version of the Birational Transformation
Theorem, extending the so-called Change of Variable Theorem for
motivic integration, due to Kontsevich \cite{kontsevich}. We have
strived to make this paper self-contained, and therefore we have
reproved the results we needed from \cite{DL}. One of our goals was
to avoid the formalism of semi-algebraic sets and work entirely in
the context of algebraic-geometry, with the hope that this will be
useful to some of the readers. In addition to the results needed for
our purpose, we have included a few other fundamental results when
we felt that our treatment simplifies the presentation available in
the literature. For example, we have included proofs of Kolchin's
Irreducibility Theorem and of Greenberg's Theorem on the
constructibility of the images of the truncation maps.

A great part of the results on spaces of jets are
characteristic--free. In particular, the Birational Transformation
Theorem holds also in positive characteristic in a form that is
slightly weaker than its usual form, but which suffices for our
applications (see Theorem~\ref{change_of_variable} below for the
precise statement). On the other hand, all our applications depend
on the existence of resolutions of singularities. Therefore we did
not shy away from using resolutions whenever this simplified the
arguments. We emphasize, however, that results such as
Theorem~\ref{theorem_introduction} above depend only on having
resolutions of singularities.

While there are no motivic integrals in these notes, the setup we discuss has
strong connections with motivic integration (in fact, the first proofs of the results connecting invariants
of singularities with spaces of arcs used this framework, see \cite{mustata}
and \cite{EMY}). For a beautiful  introduction
to the circle of ideas around motivic integration, we refer the reader to Loeser's 
Seattle lecture notes \cite{Loeser}, in this volume.

\bigskip

The paper is organized as follows. The sections \S 2--\S 6 are
devoted to the general theory of jet schemes and spaces of arcs. In
\S 2 we construct the jet schemes and prove their basic properties.
In the next section we treat the spaces of arcs and give a proof of
Kolchin's Theorem saying that in characteristic zero the space of
arcs of an irreducible variety is again irreducible. Section 4
contains two key technical results concerning the fibers of the
truncation morphisms between jet schemes. These are applied in \S 5
to study cylinders in the space of arcs of an arbitrary variety. In
particular, we prove Greenberg's Theorem and discuss the codimension
of cylinders. In \S 6 we present the Birational Transformation
Theorem of Denef-Loeser, with a simplified proof following
\cite{Lo}. This is the crucial ingredient for relating the
codimensions of cylinders in the spaces of arcs of $X'$ and of $X$,
when $X'$ is a resolution of singularities of $X$.

The reader already familiar with the basics about the codimension of
cylinders in spaces of arcs can jump directly to \S 7. Here we give
the interpretation of minimal log discrepancies from \cite{EMY}, but
without any recourse to motivic integration. In addition, we prove
our new description of these invariants in terms of contact loci in
the jet schemes. We apply this description in \S 8 to prove the
version of Inversion of Adjunction in
Theorem~\ref{theorem_introduction}. The last section is an appendix
in which we collect some general facts that we use in the main body
of the paper. In particular, in \S 9.2 we describe the connection
between the Jacobian subscheme of a variety and the subscheme
$V(J_r)$ that appears in Theorem~\ref{theorem_introduction}.

\subsection*{Acknowledgements}
The debt we owe to the paper \cite{DL} of Denef and Loeser can not
be overestimated. In addition, we have received a lot of help from
Bernd Ulrich. He explained to us the material in \S 9.2, which got
us started in our present treatment. We are grateful to Kyle Hofmann
for pointing out several typos in a preliminary version.
These notes were written
while the second author visited the Institute for Advanced Study. He
would like to thank his hosts for the stimulating environment.

\section{Jet schemes: construction and basic properties}

We work over an algebraically closed field $k$ of arbitrary
characteristic. A variety is an integral scheme, separated and of finite type over $k$. The set of nonnegative integers is denoted by $\NN$.

 Let $X$ be a scheme of finite type over $k$, and $m\in\NN$. 
 We call a scheme $J_m(X)$ over $k$ the
$m^{\rm th}$ \emph{jet scheme} of $X$ if for every $k$--algebra $A$
we have a functorial bijection
\begin{equation}\label{eq_def1}
\Hom(\Spec(A),J_m(X))\simeq\Hom(\Spec\,A[t]/(t^{m+1}),X).
\end{equation}
In particular, the $k$--valued points of $J_m(X)$ are in bijection
with the $k[t]/(t^{m+1})$--valued points of $X$. The bijections
(\ref{eq_def1}) describe the functor of points of $J_m(X)$. It
follows that if $J_m(X)$ exists, then it is unique up to  a
canonical isomorphism.

Note that if the jet schemes $J_m(X)$ and $J_p(X)$ exist and if
$m>p$, then we have a canonical projection $\pi_{m,p}\colon
J_m(X)\to J_p(X)$. This can be defined at the level of the functor
of points via (\ref{eq_def1}): the induced map
$$\Hom(\Spec\,A[t]/(t^{m+1}),X)\to\Hom(\Spec\,A[t]/(t^{p+1}),X)$$
is induced by the truncation morphism $A[t]/(t^{m+1})\to
A[t]/(t^{p+1})$. It is clear that these morphisms are compatible
whenever they are defined: $\pi_{m,p}\circ\pi_{q,m}=\pi_{q,p}$ if
$p<m<q$. If the scheme $X$ is not clear from the context, then we
write $\pi_{m,p}^X$ instead of $\pi_{m,p}$.

\begin{example}\label{example1}
We clearly have $J_0(X)=X$. For every $m$, we denote the canonical
projection $\pi_{m,0}\colon J_m(X)\to X$ by $\pi_m$.
\end{example}

\begin{proposition}\label{prop1}
For every scheme $X$ of finite type over $k$, and for every
nonnegative integer $m$, there is an $m^{\rm th}$ jet scheme
$J_m(X)$ of $X$, and this is again a scheme of finite type over $k$.
\end{proposition}

Before proving the proposition we give the following lemma.

\begin{lemma}\label{lem1}
If $U\subseteq X$ is an open subset and if $J_m(X)$ exists, then
$J_m(U)$ exists and $J_m(U)=\pi_m^{-1}(U)$.
\end{lemma}

\begin{proof}
Indeed, let $A$ be a $k$--algebra and let
$\iota_A\colon\Spec(A)\to\Spec\,A[t]/(t^{m+1})$ be induced by
truncation. Note that a morphism $f\colon \Spec\,A[t]/(t^{m+1})\to
X$ factors through $U$ if and only if the composition $f\circ
\iota_A$ factors through $U$ (factoring through $U$ is a
set-theoretic statement). Therefore the assertion of the lemma
follows from definitions.
\end{proof}

\begin{proof}[Proof of Proposition~\ref{prop1}]
Suppose first that $X$ is affine, and consider a closed embedding
$X\hookrightarrow\AAA^n$ such that $X$ is defined by the ideal
$I=(f_1,\ldots,f_q)$. For every $k$--algebra $A$, giving a morphism
$\Spec\,A[t]/(t^{m+1})\to X$ is equivalent with giving a morphism
$\phi\colon k[x_1,\ldots,x_n]/I\to A[t]/(t^{m+1})$. Such a morphism
is determined by $u_i=\phi(x_i)=\sum_{j=0}^ma_{i,j}t^j$ such that
$f_{\ell}(u_1,\ldots,u_n)=0$ for every $\ell$. We can write
$$f_{\ell}(u_1,\ldots,u_n)=\sum_{p=0}^mg_{\ell,p}((a_{i,j})_{i,j})t^p,$$
for suitable polynomials $g_{\ell,p}$ depending only on $f_{\ell}$.
It follows that $J_m(X)$ can be defined in $\AAA^{(m+1)n}$ by the
polynomials $g_{\ell,p}$ for $1\leq \ell\leq q$ and $0\leq p\leq m$.

Suppose now that $X$ is an arbitrary scheme of finite type over $k$.
Consider an affine cover $X=U_1\cup\ldots\cup U_r$. As we have seen,
we have an $m^{\rm th}$ jet scheme $\pi_m^i\colon J_m(U_i)\to U_i$
for every $i$. Moreover, by Lemma~\ref{lem1}, for every $i$ and $j$,
the inverse images $(\pi_m^i)^{-1}(U_i\cap U_j)$ and
$(\pi_m^j)^{-1}(U_i\cap U_j)$ give the $m^{\rm th}$ jet scheme of
$U_i\cap U_j$. Therefore they are canonically isomorphic. This shows
that we may construct a scheme $J_m(X)$ by glueing the schemes
$J_m(U_i)$ along the canonical isomorphisms of
$(\pi_m^i)^{-1}(U_i\cap U_j)$ with $(\pi_m^j)^{-1}(U_i\cap U_j)$.
Moreover, the projections $\pi_m^i$ also glue to give a morphism
$\pi_m\colon J_m(X)\to X$. It is now straightforward to check that
$J_m(X)$ has the required property.
\end{proof}

\begin{remark}
It follows from the description in the above proof that for every
$X$, the projection $\pi_m\colon J_m(X)\to X$ is affine.
\end{remark}

\begin{example}\label{example2}
The first jet scheme $J_1(X)$ is isomorphic to the total tangent
space $TX:={\mathcal Spec}({\rm Sym}(\Omega_{X/k}))$. Indeed, arguing as
in the proof of Proposition~\ref{prop1}, we see that it is enough to
show the assertion when $X=\Spec(R)$ is affine, in which case
$TX={\rm Spec}({\rm Sym}(\Omega_{R/k}))$. In this case, if $A$
is a $k$--algebra, then giving a  morphism of schemes
$f\colon\Spec(A)\to \Spec({\rm Sym}(\Omega_{R/k}))$ is equivalent
with giving a morphism of $k$--algebras $\phi\colon R\to A$ and a
$k$-derivation $D\colon R\to A$ (where $A$ becomes an $R$-module via
$\phi$). This is the same as giving a ring homomorphism $f\colon R
\to A[t]/(t^2)$, where $f(u)=\phi(u)+tD(u)$.
\end{example}

If $f\colon X\to Y$ is a morphism of schemes, then we get a
corresponding morphism $f_m\colon J_m(X)\to J_m(Y)$. At  the level
of $A$--valued points, this takes an $A[t]/(t^{m+1})$--valued point
$\gamma$ of $X$ to $f\circ\gamma$. Taking $X$ to $J_m(X)$ gives a
functor from the category of schemes of finite type over $k$ to
itself. Note also that the morphisms $f_m$ are compatible in the
obvious sense with the projections $J_m(X)\to J_{m-1}(X)$ and
$J_m(Y)\to J_{m-1}(Y)$.

\begin{remark}\label{affine_space}
The jet schemes of the affine space are easy to describe: we have an
isomorphism $J_m(\AAA^n)\simeq\AAA^{(m+1)n}$ such that the
projection $J_m(\AAA^n)\to J_{m-1}(\AAA^n)$ corresponds
to the projection onto the first $mn$ coordinates. Indeed, an
$A$--valued point of $J_m(\AAA^n)$ corresponds to a ring
homomorphism $\phi\colon k[x_1,\ldots,x_n]\to A[t]/(t^{m+1})$, which
is uniquely determined by giving each $\phi(X_i)\in
A[t]/(t^{m+1})\simeq A^{m+1}$.
\end{remark}

\begin{remark}\label{closed_embedding}
In light of the previous remark, we see that the proof of
Proposition~\ref{prop1} showed that if $i\colon X\hookrightarrow
\AAA^n$ is a closed immersion, then the induced morphism $i_m\colon
J_m(X)\to J_m(\AAA^n)$ is also a closed immersion. Using the
description of the equations of $J_m(X)$ in $J_m(\AAA^n)$ we see
that more generally, if $f\colon X\hookrightarrow Y$ is a closed
immersion, then $f_m$ is a closed immersion, too.
\end{remark}

\begin{remark}
The following are some direct consequences of the definition.
\begin{enumerate}
\item[i)] For every schemes $X$ and $Y$ and for every $m$, there is a
canonical isomorphism $J_m(X\times Y)\simeq J_m(X)\times J_m(Y)$.
\item[ii)] If $G$ is a group scheme over $k$, then $J_m(G)$ is also a group
scheme over $k$. Moreover, if $G$ acts on $X$, then $J_m(G)$ acts on
$J_m(X)$.
\item[iii)] If $f\colon Y\to X$ is a morphism of schemes and $Z\hookrightarrow X$
is a closed subscheme, then we have a canonical isomorphism
$J_m(f^{-1}(Z))\simeq f_m^{-1}(J_m(Z))$.
\end{enumerate}
\end{remark}

The following lemma generalizes Lemma~\ref{lem1} to the case of an
\'{e}tale morphism.

\begin{lemma}\label{lem2}
If $f\colon X\to Y$ is an \'{e}tale morphism, then for every $m$ the
commutative diagram

\[
\begin{CD}
J_m(X) @>{f_m}>> J_m(Y) \\
@VV{\pi_m^X}V @VV{\pi_m^Y}V \\
X @>{f}>>Y
\end{CD}
\]
is Cartesian.
\end{lemma}

\begin{proof}
{}From the description of the $A$--valued points of $J_m(X)$ and
$J_m(Y)$ we see that it is enough to show  that for every
$k$--algebra $A$ and every commutative diagram
\[
\begin{CD}
\Spec(A)@>>> X\\
@VVV@VVV\\
\Spec\,A[t]/(t^{m+1}) @>>>Y
\end{CD}
\]
there is a unique morphism $\Spec\,A[t]/(t^{m+1})\to X$ making the
two triangles commutative. This is a consequence of the fact that
$f$ is formally \'{e}tale.
\end{proof}

\begin{remark}
A similar argument shows that if $f\colon Y\to X$ is a smooth
surjective morphism, then $f_m$ is surjective for every $m$.
Moreover, $f_m$ is again smooth: this follows from Lemma~\ref{lem2}
and the fact that $f$ can be locally factored as $U\overset{g}\to
V\times{\mathbb A}^n\overset{p}\to V$, where $g$ is \'{e}tale and
$p$ is the projection onto the first component.
\end{remark}

We say that a morphism of schemes $g\colon V'\to V$ is \emph{locally
trivial} with fiber $F$ if there is a cover by Zariski open subsets
$V=U_1\cup\ldots\cup U_r$ such that $g^{-1}(U_i)\simeq U_i\times F$,
with the restriction of $g$ corresponding to the projection onto the
first component.

\begin{corollary}\label{cor1}
If $X$ is a nonsingular variety of dimension $n$, then all
projections $\pi_{m,m-1}\colon J_m(X)\to J_{m-1}(X)$ are locally
trivial with fiber $\AAA^n$. In particular, $J_m(X)$ is a nonsingular
variety of dimension $(m+1)n$.
\end{corollary}

\begin{proof}
Around every point in $X$ we can find an open subset $U$ and an
\'{e}tale morphism $U\to\AAA^n$. Using Lemma~\ref{lem2} we reduce
our assertion to the case of the affine space, when it follows from
Remark~\ref{affine_space}.
\end{proof}

\begin{remark}
If $X$ and $Y$ are schemes and $x\in X$ and $y\in Y$ are points such
that the completions $\widehat{\OO}_{X,x}$ and $\widehat{\OO}_{Y,y}$
are isomorphic, then the fiber of $J_m(X)$ over $x$ is isomorphic to
the fiber of $J_m(Y)$ over $y$. Indeed, the $A$--valued points of
the fiber of $J_m(X)$ over $x$ are in natural bijection with
$$\{\phi\colon\OO_{X,x}\to A[t]/(t^{m+1})\mid\phi(\frmm_x)\subseteq
(t)\}=\{\hat{\phi}\colon\widehat{\OO}_{X,x}\to A[t]/(t^{m+1})\mid
\hat{\phi}(\widehat{\frmm}_x)\subseteq (t)\}$$
$$\simeq\{\hat{\psi}\colon\widehat{\OO}_{Y,y}\to A[t]/(t^{m+1})\mid
\hat{\psi}(\widehat{\frmm}_y)\subseteq
(t)\}=\{\psi\colon\OO_{Y,y}\to
A[t]/(t^{m+1})\mid\psi(\frmm_y)\subseteq(t)\}.$$
\end{remark}

\begin{example}
Suppose that $X$ is a reduced curve having a node at $p$, i.e. we
have ${\widehat{\mathcal O}}_{X,p}\simeq k\llbracket x,y\rrbracket/(xy)$. By the
previous remark, in order to compute the fiber of $J_m(X)$ over $p$
we may assume that $X=\Spec\,k[x,y]/(xy)$ and that $p$ is the
origin.
 We see that this fiber consists of the union of $m$ irreducible
 components, each of them (with the reduced structure) being isomorphic to
 $\AAA^{m+1}$. Indeed, the $i^{\rm th}$ such component corresponds to
 morphisms $\phi\colon k[x,y]\to k[t]/(t^{m+1})$ such that
 ${\rm ord}(\phi(x))\geq i$ and ${\rm ord}(\phi(y))\geq m+1-i$.

 If $C$ is an irreducible component of $X$ passing
 through $p$ and $C_{\rm reg}$ is its nonsingular locus, then Corollary~\ref{cor1} implies that
 $\overline{J_m(C_{\rm reg})}$ is an irreducible component of
 $J_m(X)$ of dimension $(m+1)$. Therefore all the above components
 of the fiber of $J_m(X)$ over $p$ are irreducible components of
 $J_m(X)$. In particular, $J_m(X)$ is not irreducible for every
 $m\geq 1$.
\end{example}

\begin{example}
Let $X$ be an arbitrary scheme and $p$ a point in $X$. If all
projections $(\pi^X_m)^{-1}(p)\to (\pi^X_{m-1})^{-1}(p)$ are
surjective, then $p$ is a nonsingular point. To see this, it is
enough to show that if a tangent vector in $T_pX$ can be lifted to
any $J_m(X)$, then it lies in the tangent cone of $X$ at $p$. We may
assume that $X$ is a closed subscheme of $\AAA^n$ and that $p$ is
the origin. The tangent cone of $X$ at $p$ is the intersection of
the tangent cone at $p$ to each hypersurface $H$ containing $X$.
Since $J_m(X)\subseteq J_m(H)$ for every $m$ and every such $H$, it
is enough to prove our assertion when $X$ is a hypersurface. Let $f$
be an equation defining $X$, and write $f=f_r+f_{r+1}+\ldots$, where
$f_i$ has degree $i$ and $f_r\neq 0$. By considering the equations
defining $J_r(X)$ in $J_r(\AAA^n)$, we see that the commutative diagram
\[
\begin{CD}
(\pi^X_r)^{-1}(p) @>>> (\pi_r^{\AAA^n})^{-1}(p)=\AAA^n\times\AAA^{(r-1)n} \\
@VVV @VV{\rm pr}_1V \\
T_pX @>>>T_p\AAA^n=\AAA^n
\end{CD}
\]
identifies the fiber of $J_r(X)$ over $p$ with
$T\times\AAA^{(r-1)n}\hookrightarrow \AAA^n\times\AAA^{(r-1)n}$,
where $T$ is defined by $f_r$
in $\AAA^n$. Since $T$ is the tangent cone to $X$ at $p$, this
completes the proof of our assertion.
\end{example}

\section{Spaces of arcs}

We now consider the projective limit of the jet schemes. Suppose
that $X$ is a scheme of finite type over $k$. Since the projective
system
$$\cdots\to J_m(X)\to J_{m-1}(X)\to\cdots\to J_0(X)=X$$
consists of affine morphisms, the projective limit exists in the
category of schemes over $k$. It is denoted by $J_{\infty}(X)$ and
it is called the space of arcs of $X$. In general, it is not of
finite type over $k$.

The space of arcs comes equipped with projection morphisms
$\psi_m\colon J_{\infty}(X)\to J_m(X)$ that are affine. In
particular, we have $\psi_0\colon J_{\infty}(X)\to X$. Over an
affine open subset $U\subseteq X$, the space of arcs is described by
$$\OO(\psi_0^{-1}(U))=\underrightarrow{\rm lim}\,\OO(\pi_m^{-1}(U)).$$

It follows from the projective limit definition and the functorial
description of the jet schemes that if $X$ is affine, then for every
$k$--algebra $A$ we have
\begin{equation}
\Hom(\Spec(A), J_{\infty}(X))\simeq
\underleftarrow\Hom(\Spec\,A[t]/(t^{m+1}), X)
\simeq\Hom(\Spec\,A\llbracket t\rrbracket,X).
\end{equation}
If $X$ is not necessarily affine, note that every morphism
$\Spec\,k[t]/(t^{m+1})\to X$ or $\Spec\,k\llbracket t\rrbracket\to X$ factors through
any affine open neighborhood of the image of the closed point. It
follows that for every $X$, the $k$--valued points of
$J_{\infty}(X)$ correspond to \emph{arcs in $X$}
$$\Hom(\Spec(k),J_{\infty}(X))\simeq\Hom(\Spec\,k\llbracket t\rrbracket, X).$$

If $f\colon X\to Y$ is a morphism of schemes, by taking the
projective limit of the morphisms $f_m$ we get a morphism
$f_{\infty}\colon J_{\infty}(X)\to J_{\infty}(Y)$. We get in this
way a functor from $k$-schemes of finite type over $k$ to arbitrary
$k$-schemes (in fact, to quasicompact and quasiseparated
$k$-schemes).

The properties we have discussed in the previous section for jet
schemes induce corresponding properties for the spaces of arcs. For
example, if $f\colon X\to Y$ is an \'{e}tale morphism, then we have
a Cartesian diagram

\[
\begin{CD}
J_{\infty}(X) @>{f_{\infty}}>> J_{\infty}(Y) \\
@VV{\psi_0^X}V @VV{\psi_0^Y}V \\
X @>{f}>>Y.
\end{CD}
\]

If $i\colon X\hookrightarrow Y$ is a closed immersion, then
$i_{\infty}$ is also a closed immersion. Moreover, if $Y=\AAA^n$,
then $J_{\infty}(Y)\simeq\AAA^{\NN}=\Spec\,k[x_1,x_2,\ldots]$, such that $\psi_m$
corresponds to the projection onto the first $(m+1)n$ components. As
in the proof of Proposition~\ref{prop1}, starting with equations for
a closed subscheme $X$ of $\AAA^n$ we can write down equations for
$J_{\infty}(X)$ in $J_{\infty}(\AAA^n)$.

\smallskip

Note that the one-dimensional torus $k^*$ has a natural action on
jet schemes induced by reparametrization the jets. In fact, for
every scheme $X$ we have a morphism
$$\Phi_m\colon\AAA^1\times J_m(X)\to J_m(X)$$
described at the level of functors of points as follows. For every
$k$--algebra $A$, an $A$--valued point of $\AAA^1\times J_m(X)$
corresponds to a pair $(a,\phi)$, where $a\in A$ and $\phi\colon
{\rm Spec}\,A[t]/(t^{m+1})\to X$. This pair is mapped by $\Phi_m$ to
the $A$--valued point of $J_m(X)$ given by the composition
$$\Spec\,A[t]/(t^{m+1})\to\Spec\,A[t]/(t^{m+1})\overset{\phi}\to
X,$$ where the first arrow corresponds to the ring homomorphism
induced by $t\to at$.

It is clear that $\Phi_m$ induces an action of $k^*$ on $J_m(X)$.
The fixed points of this action are given by $\Phi_m(\{0\}\times
J_m(X))$. These are the \emph{constant jets} over the points in $X$:
over a point $x\in X$ the constant $m$--jet is the composition
$$\gamma_m^x\colon \Spec\,k[t]/(t^{m+1})\to\Spec\,k\to X,$$
where the second arrow gives $x$. We have a \emph{zero-section}
$s_m\colon X\to J_m(X)$ of the projection $\pi_m$ that takes $x$ to
$\gamma_m^x$. If $A$ is a $k$--algebra, then $s_m$ takes an
$A$--valued point of $X$ given by $u\colon {\rm Spec}\,A\to X$ to
the composition
$${\rm Spec}\,A[t]/(t^{m+1})\to {\rm Spec}\,A\overset{u}\to X,$$
the first arrow being induced by the inclusion $A\hookrightarrow
A[t]/(t^{m+1})$.

Note that if $\gamma\in J_m(X)$ is a jet lying over $x\in X$, then
$\gamma_m^x$ lies in the closure of $\Phi_m(k^*\times\{\gamma\})$.
Since every irreducible component $Z$ of $J_m(X)$ is preserved by
the $k^*$--action, this implies that if $\gamma$ is an $m$--jet in
$Z$ that lies over $x\in X$, then also $\gamma_m^x$ is in $Z$. This
will be very useful for the applications in \S 8.

Both the morphisms $\Phi_m$ and the zero-sections $s_m$ are
functorial. Moreover, they satisfy obvious compatibilities with the
projections $J_m(X)\to J_{m-1}(X)$. Therefore we get a morphism
$$\Phi_{\infty}\colon\AAA^1\times J_{\infty}(X)\to J_{\infty}(X)$$
inducing an action of $k^*$ on $J_{\infty}(X)$, and a zero-section
$s_{\infty}\colon X\to J_{\infty}(X)$.

\bigskip

If ${\rm char}(k)=0$, then one can write explicit equations
 for $J_{\infty}(X)$ and $J_m(X)$ by "formally differentiating", as follows. If $S=k[x_1,\ldots,x_n]$,
let us write $S_{\infty}=k[x_i^{(m)}\mid 1\leq i\leq n, m\in\NN]$,
so that $\Spec(S_{\infty})=J_{\infty}(\AAA^n)$ (in practice, we
simply write $x_i=x_i^{(0)}$, $x_i'=x_i^{(1)}$, and so on). 
The identification is made as follows: for a
$k$--algebra $A$, a morphism $\phi\colon
k[x_1,\ldots,x_n]\to A\llbracket t\rrbracket$ determined by
\begin{equation}\label{def_phi}
\phi(x_i)=\sum_{m\in\NN}\frac{a_i^{(m)}}{m!}t^m
\end{equation}
corresponds to the $A$--valued point $(a_i^{(m)})$ of $\Spec(S_{\infty})$.

Note that on $S_{\infty}$ we have a $k$--derivation $D$ characterized by
$D(x_i^{(m)})=x_i^{(m+1)}$.
If $f\in R$, then we put $f^{(0)}:=f$, and we define recursively
$f^{(m)}:=D(f^{(m-1)})$ for $m\geq 1$. Suppose now that $R=S/I$, where $I$ is
generated by $f_1,\ldots,f_r$. We claim that if
\begin{equation}
R_{\infty}:=S_{\infty}/(f_i^{(m)}\vert 1\leq i\leq r, m\in\NN),
\end{equation}
then $J_{\infty}(\Spec\,R)\simeq \Spec(R_{\infty})$.

Indeed, given $A$ and $\phi$ as above, for every
$f\in k[x_1,\ldots,x_n]$ we have
$$\phi(f)=\sum_{m\in\NN}\frac{f^{(m)}(a,a',\ldots,a^{(m)})}{m!}t^m$$
(note that both sides are additive and multiplicative in $f$, hence
it is enough to check this for $f=x_i$, when it is trivial). It
follows that $\phi$ induces a morphism $R\to A\llbracket t\rrbracket$ if and only if
$f_i^{(m)}(a,a',\ldots,a^{(m)})=0$ for every $m$ and every $i\leq
r$. This completes the proof of the above claim.

\begin{remark}
Note that $D$ induces a $k$--derivation $\overline{D}$ on $R_{\infty}$.
Moreover, $(R_{\infty},\overline{D})$ is universal in the following
sense: we have a $k$--algebra homomorphism $j\colon R\to R_{\infty}$
such that if $(T,\delta)$ is another $k$--algebra with a
$k$--derivation $\delta$, and if $j'\colon R\to T$ is a $k$--algebra
homomorphism, then there is a unique $k$--algebra homomorphism
$h\colon R_{\infty}\to T$ making the diagram
\[
\xymatrix{
R\ar[dr] _{j'}     & \overset{j}\longrightarrow& (R_{\infty},\overline{D}) \ar[dl]^{h}  \\
& (T,\delta)
}
\]
commutative, and such that $h$
commutes with the derivations, i.e.
$\delta(h(u))=h(\overline{D}(u))$ for every $u\in R_{\infty}$. This
is the starting point for the applications of spaces of arcs in
differential algebra, see \cite{buium}.
\end{remark}

Of course, if we consider finite level truncations, then we obtain
equations for the jet schemes. More precisely, if we put
$S_m:=k[x_i^{(j)}\mid i\leq n, 0\leq j\leq m]$ and
$$R_m:=S_m/(f_i,f_i',\ldots,f_i^{(m)}\mid 1\leq i\leq r),$$
then $\Spec(R_m)\simeq J_m(\Spec\,R)$. Moreover, the 
obvious morphisms $R_{m-1}\to R_m$
induce the projections
$J_m(\Spec\,R)\to J_{m-1}(\Spec\,R)$.

\smallskip

{}From now on, whenever dealing with the schemes $J_m(X)$ and
$J_{\infty}(X)$ we will restrict to their $k$--valued points. Of
course, for $J_m(X)$ this causes no ambiguity since this is a scheme
of finite type over $k$. Note that the Zariski topology on
$J_{\infty}(X)$ is the projective limit topology of
$J_{\infty}(X)\simeq\underleftarrow{\rm lim} J_m(X)$. Moreover,
since we consider only $k$--valued points, we have
$J_{\infty}(X)=J_{\infty}(X_{\rm red})$ (note that the analogous
assertion is false for the spaces $J_m(X)$). Indeed, since $k\llbracket t\rrbracket$
is a domain, we have ${\rm Hom}(\Spec\,k\llbracket t\rrbracket,X)={\rm
Hom}(\Spec\,k\llbracket t\rrbracket, X_{\rm red})$. Similarly, if
$X=X_1\cup\ldots\cup X_r$, where all $X_i$ are closed in $X$, then
$J_{\infty}(X)=J_{\infty}(X_1)\cup\ldots \cup J_{\infty}(X_r)$.

\bigskip

To a closed subscheme $Z$  of a scheme $X$ we associate subsets of
the spaces of arcs and jets of $X$ by specifying the vanishing order
along $Z$. If $\gamma\colon\Spec\,k\llbracket t\rrbracket \to X$ is an arc on $X$,
then the inverse image of $Z$ by $\gamma$ is defined by an ideal in
$k\llbracket t\rrbracket$. If this ideal is generated by $t^r$, then we put ${\rm
ord}_{\gamma}(Z)=r$ (if the ideal is zero, then we put ${\rm
ord}_{\gamma}(Z)=\infty$). The \emph{contact locus} of order $e$
with $Z$ in $J_{\infty}(X)$ is the set
$${\rm Cont}^e(Z):=\{\gamma\in
J_{\infty}(X)\mid\ord_{\gamma}(Z)=e\}.$$ We similarly define
$${\rm Cont}^{\geq e}(Z):=\{\gamma\in
J_{\infty}(X)\mid\ord_{\gamma}(Z)\geq e\}.$$ We can define in the
obvious way also subsets ${\rm Cont}^e(Z)_m$ (if $e\leq m$) and
${\rm Cont}^{\geq e}(Z)_m$ (if $e\leq m+1$) of $J_m(X)$ and we have
$${\rm Cont}^e(Z)=\psi_m^{-1}({\rm
Cont}^e(Z)_m),\,\,{\rm Cont}^{\geq e}(Z)=\psi_m^{-1}({\rm
Cont}^{\geq e}(Z)_m).$$ Note that ${\rm Cont}^{\geq (m+1)}(Z)_m
=J_m(Z)$. If ${\mathcal I}$ is the ideal sheaf in $\OO_X$ defining
$Z$, then we sometimes write ${\rm ord}_{\gamma}({\mathcal I})$,
${\rm Cont}^e({\mathcal I})$ and ${\rm Cont}^{\geq e}({\mathcal
I})$.

\bigskip

The next proposition gives the first hint of the relevance of spaces
of arcs to birational geometry. A key idea is that certain subsets
in the space of arcs are "small" and they can be ignored. A subset
of $J_{\infty}(X)$ is called \emph{thin} if it is contained in
$J_{\infty}(Y)$, where $Y$ is a closed subset of $X$ that does not
contain an irreducible component of $X$. It is clear that a finite
union of thin subsets is again thin. If $f\colon X'\to X$ is a
dominant morphism with $X$ and $X'$ irreducible, and $A\subseteq
J_{\infty}(X)$ is thin, then $f_{\infty}^{-1}(A)$ is thin. If $f$ is
in addition generically finite, and $B\subseteq J_{\infty}(X')$ is
thin, then $f_{\infty}(B)$ is thin.

We show that a proper birational morphism induces a bijective map on
the complements of suitable thin sets.

\begin{proposition}\label{prop11}
Let $f\colon X'\to X$ be a proper morphism. If $Z$ is a closed
subset of $X$ such that $f$ is an isomorphism over $X\smallsetminus
Z$, then the induced map
$$J_{\infty}(X')\smallsetminus J_{\infty}(f^{-1}(Z))\to
J_{\infty}(X)\smallsetminus J_{\infty}(Z)$$ is bijective. In
particular, if $f$ is a proper birational morphism of reduced
schemes, then $f_{\infty}$ gives a bijection on the complements of
suitable thin subsets.
\end{proposition}

\begin{proof}
Let $U=X\smallsetminus Z$. Since $f$ is proper, the Valuative
Criterion for Properness implies that an arc $\gamma\colon
\Spec\,k\llbracket t\rrbracket \to X$ lies in the image of $f_{\infty}$ if and only if
the induced morphism $\overline{\gamma}\colon\Spec\,k\llparenthesis t\rrparenthesis\to X$ can
be lifted to $X'$ (moreover, if the lifting of $\overline{\gamma}$
is unique, then the lifting of $\gamma$ is also unique). On the
other hand, $\gamma$ does not lie in $J_{\infty}(Z)$ if and only if
$\overline{\gamma}$ factors through $U\hookrightarrow X$. In this
case, the lifting of $\overline{\gamma}$ exists and is unique since
$f$ is an isomorphism over $U$.
\end{proof}

We use the above proposition to prove the following result of
Kolchin.

\begin{theorem}\label{thm1}
If $X$ is irreducible and ${\rm char}(k)=0$, then $J_{\infty}(X)$ is
irreducible.
\end{theorem}

\begin{proof}
Since $J_{\infty}(X)=J_{\infty}(X_{\rm red})$, we may assume that
$X$ is also reduced.
 If $X$ is nonsingular, then the assertion in the theorem is easy: we have
seen that every jet scheme $J_m(X)$ is a nonsingular variety. Since
the projections $J_{\infty}(X)\to J_m(X)$ are surjective, and
$J_{\infty}(X)=\underleftarrow{\lim}J_m(X)$ with the projective
limit topology, it follows that $J_{\infty}(X)$, too,  is
irreducible.

In the general case we do induction on $n=\dim(X)$, the case $n=0$
being trivial. By Hironaka's Theorem we have a resolution of
singularities $f\colon X'\to X$, that is, a proper birational
morphism, with $X'$ nonsingular. Suppose that $Z$ is a proper closed
subset of $X$ such that $f$ is an isomorphism over
$U=X\smallsetminus Z$. It follows from Proposition~\ref{prop11} that
$$J_{\infty}(X)=J_{\infty}(Z)\cup {\rm Im}(f_{\infty}).$$
Moreover, the nonsingular case implies that $J_{\infty}(X')$, hence
also ${\rm Im}(f_{\infty})$, is irreducible. Therefore, in order to
complete the proof it is enough to show that $J_{\infty}(Z)$ is
contained in the closure of ${\rm Im}(f_{\infty})$.

Consider the irreducible decomposition $Z=Z_1\cup\ldots\cup Z_r$,
inducing $J_{\infty}(Z)=J_{\infty}(Z_1)\cup\ldots \cup
J_{\infty}(Z_r)$. Since $f$ is surjective, for every $i$ there is an
irreducible component $Z'_i$ of $f^{-1}(Z_i)$ such that the induced
map $Z'_i\to Z_i$ is surjective. We are in characteristic zero,
hence by the Generic Smoothness Theorem we can find open subsets
$U'_i$ and $U_i$ in $Z'_i$ and $Z_i$, respectively, such that the
induced morphisms $g_i\colon U'_i\to U_i$ are smooth and surjective.
In particular, we have
$$J_{\infty}(U_i)={\rm Im}((g_i)_{\infty})\subseteq {\rm
Im}(f_{\infty}).$$

On the other hand, every $J_{\infty}(Z_i)$ is irreducible by
induction. Since $J_{\infty}(U_i)$ is a nonempty open subset of
$J_{\infty}(Z_i)$, it follows that
$$J_{\infty}(Z_i)\subseteq \overline{{\rm Im}(f_{\infty})}$$
for every $i$. This completes the proof of the theorem.
\end{proof}

\begin{remark}
In fact, Kolchin's Theorem holds in a much more general setup, see
\cite{kolchin} and also \cite{gillet} for a scheme-theoretic
aproach. In fact, we proved a slightly weaker statement even in our
restricted setting. Kolchin's result says that the scheme
$J_{\infty}(X)$ is irreducible, while we only proved  that its
$k$--valued points form an irreducible set. In fact, one can deduce
the stronger statement from ours by showing that the $k$--valued
points are dense in $J_{\infty}(X)$. In turn, this can be proved in
a similar way with Theorem~\ref{thm1} above. For a different proof
of (the stronger version of) Kolchin's Theorem, without using
resolution of singularities, see \cite{IK} and \cite{NS}. Note also
that Remark~1 in \cite{NS} gives a counterexample in positive
characteristic.
\end{remark}

\section{Truncation maps between spaces of jets}

In what follows we will encounter morphisms that are not locally
trivial, but that satisfy this property after passing to a
stratification. Suppose that $g\colon V'\to V$ is a morphism of
schemes, $W'\subseteq V'$ and $W\subseteq V$ are constructible
subsets such that $g(W')\subseteq W$, and $F$ is a reduced scheme.
We will say that $g$ gives a \emph{piecewise trivial} fibration
$W'\to W$ with fiber $F$ if there is a decomposition
$W=T_1\sqcup\ldots\sqcup T_r$, with all $T_i$ locally closed
subsets of $W$ (with the reduced scheme structure) such that each
$W'\cap g^{-1}(T_i)$ is locally closed in $V$ and, with the reduced
scheme structure, it is isomorphic to $T_i\times F$ (with the
restriction of $g$ corresponding to the projection onto the first
component). It is clear that if $g\colon V'\to V$ is locally trivial
with fiber $F$, then it gives a piecewise trivial fibration with
fiber $F_{\rm red}$ from $g^{-1}(W)$ to $W$ for every constructible
subset $W$ of $V$.

If in the definition of piecewise trivial fibrations we assume only
that $W'\cap g^{-1}(T_i)\to T_i$ factors as
$$W'\cap g^{-1}(T_i)\overset{u}\to T'_i\times F\overset{v}\to T'_i\overset{w}\to T_i,$$
where $u$ is an isomorphism, $v$ is the projection, and $w$ is
bijective, then we say that $W'\to W$ is a \emph{weakly piecewise
trivial} fibration with fiber $F$. If ${\rm char}(k)=0$, then every
bijective morphism is piecewise trivial with fiber ${\rm Spec}(k)$,
and therefore the two notions coincide.

\smallskip

We have seen in Corollary~\ref{cor1} that if $X$ is a nonsingular
variety of dimension $n$, then the truncation maps $J_m(X)\to
J_{m-1}(X)$ are locally trivial with fiber $\AAA^n$. In order to
generalize this to more general schemes, we need to introduce the
\emph{Jacobian subscheme}. If $X$ is a scheme of pure dimension $n$,
then its Jacobian subscheme is defined by ${\rm Jac}_X$, the Fitting
ideal ${\rm Fitt}^n(\Omega_X)$.  For the basics on Fitting ideals we
refer to \cite{eisenbud}. A basic property of Fitting ideals that we
will keep using is that they commute with pull-back: if $f\colon
X'\to X$ is a morphism and if ${\mathcal M}$ is a coherent sheaf on
$X$, then ${\rm Fitt}^i(f^*{\mathcal M})=({\rm
Fitt}^i({\mathcal M}))\cdot {\mathcal O}_{X'}$ for every $i$.

The ideal ${\rm Jac}_X$ can be explicitly computed as follows.
Suppose that $U$ is an open subset of $X$ that admits a closed
immersion $U\hookrightarrow \AAA^N$. We have a surjection
$$\Omega_{\AAA^N}\vert_X=\oplus_{j=1}^N\OO_Xdx_j\to\Omega_X$$
with the kernel generated by the $df=\sum_{j=1}^N\frac{\partial
f}{\partial x_j}dx_j$, where $f$ varies over a system of generators
$f_1,\ldots,f_d$ for the ideal of $U$ in $\AAA^N$. If $r=N-n$, then
${\rm Jac}_X$ is generated over $U$ by the image in $\OO_U$ of the
$r$--minors of the Jacobian matrix $(\partial f_i/\partial
x_j)_{i,j}$.

 It is well-known that the support of the Jacobian subscheme is the
singular locus $X_{\rm sing}$ of $X$. Most of the time we will
assume that $X$ is reduced, hence its singular locus does not
contain any irreducible component of $X$. Note also that
${\rm Fitt}^{n-1}(\Omega_X)=0$ if either
$X$ is locally a complete intersection (when the above Jacobian matrix
has $r$ rows) or if $X$ is reduced (when the $(r+1)$--minors of the
Jacobian matrix vanish at the generic points of the irreducible components of $X$,
hence are zero in $\OO_X$).

\smallskip

We start by describing the fibers of the truncation morphisms when
we restrict to jets that can be lifted to the space of arcs.

\begin{proposition}\label{fiber2}{\rm (}\cite{DL}{\rm )}
Let $X$ be a reduced scheme of pure dimension $n$ and $e$ a
nonnegative integer. Fix $m\geq e$ and let $\pi_{m+e,m}\colon
J_{m+e}(X)\to J_m(X)$ be the canonical projection.
\begin{enumerate}
\item[i)] We have $\psi_m({\rm Cont}^e({\rm
Jac}_X))=\pi_{m+e,m}({\rm Cont}^e({\rm Jac}_X)_{m+e})$, i.e. an
$m$--jet on $J_m(X)$ that vanishes with order $e$ along ${\rm
Jac}_X$ can be lifted to $J_{\infty}(X)$ if and only if it can be
lifted to $J_{m+e}(X)$. In particular, $\psi_m({\rm Cont}^e({\rm
Jac}_X))$ is a constructible set.
\item[ii)] The projection $J_{m+1}(X)
\to J_m(X)$ induces a piecewise trivial fibration
$$\alpha\colon \psi_{m+1}({\rm Cont}^e({\rm Jac}_X))\to
\psi_m({\rm Cont}^e({\rm Jac}_X))$$ with fiber $\AAA^n$.
\end{enumerate}
\end{proposition}

Before giving the proof of Proposition~\ref{fiber2} we make some
general considerations that will be used again later. A key point
for the proof of Proposition~\ref{fiber2} is the reduction to the
complete intersection case. We present now the basic setup, leaving
the proof of a technical result for the Appendix.

Let $X$ be a reduced scheme of pure dimension $n$. All our
statements are local over $X$, hence we may assume that $X$ is
affine. Fix a closed embedding $X\hookrightarrow\AAA^N$ and let
$f_1,\ldots,f_d$ be generators of the ideal $I_X$ of $X$. Consider
$F_1,\ldots,F_d$ with $F_i=\sum_{j=1}^da_{i,j}f_j$ for general
$a_{i,j}\in k$. Note that we still have $I_X=(F_1,\ldots,F_d)$, but
in addition we have the following properties. Let us denote by $M$
the subscheme defined by the ideal $I_M=(F_1,\ldots,F_r)$, where
$r=N-n$.

\begin{enumerate}
\item[1)] All irreducible components of $M$ have dimension $n$,
hence $M$ is a complete intersection.
\item[2)] $X$ is a closed subscheme of $M$ and $X=M$ at the generic
point of every irreducible component of $X$.
\item[3)] There is an $r$--minor of the Jacobian matrix of
$F_1,\ldots,F_r$ that does not vanish at the generic point of any
irreducible component of $X$.
\end{enumerate}
Of course, every $r$ elements of $\{F_1,\ldots,F_d\}$ satisfy
analogous properties.

Suppose now that $e$ is a nonnegative integer, $m\geq e$ and we want
to study ${\rm Cont}^e({\rm Jac_X})_m$. If $M$ is as above, then we
have an open subset $U_M$ of ${\rm Cont}^e({\rm Jac}_X)_m$ that is
contained in ${\rm Cont}^e({\rm Jac}_M)_m$ (the latter contact locus
is a subset of $J_m(M)$). Moreover, when varying the subsets of
$\{1,\ldots,d\}$ with $r$ elements, the corresponding open subsets
cover ${\rm Cont}^e({\rm Jac}_X)_m$.

\begin{lemma}\label{lemma_reduction}
If $\gamma\in {\rm Cont}^e({\rm Jac}_M)\subseteq J_{\infty}(M)$ is
such that its projection to $J_m(M)$ lies in $J_m(X)$, then $\gamma$
lies in $J_{\infty}(X)$.
\end{lemma}

\begin{proof}
Let $X'\subseteq\AAA^N$ be defined by $(I_M\colon I_X)$, hence
set-theoretically $X'$ is the union of the irreducible components in
$M$ that are not contained in $X$. We have $J_{\infty}(M)=
J_{\infty}(X)\cup J_{\infty}(X')$, and therefore it is enough to
show that $\gamma$ does not lie in $J_{\infty}(X')$.

It follows from Corollary~\ref{cor1_appendix} in the Appendix that
if we denote by $J_F$ the ideal generated by the $r$--minors of the
Jacobian matrix of $(F_1,\ldots,F_r)$ (hence ${\rm
Jac}_M=(J_F+I_M)/I_M$), then
$$J_F\subseteq I_{X'}+I_X.$$
By assumption ${\rm ord}_{\gamma}(J_F)=e<m+1\leq {\rm
ord}_{\gamma}(I_X)$, hence ${\rm ord}_{\gamma}(I_{X'})\leq e$. In
particular, $\gamma$ is not in $J_{\infty}(X')$.
\end{proof}

\smallskip

\begin{proof}[Proof of Proposition~\ref{fiber2}]

We may assume that $X$ is affine, and let $X\hookrightarrow \AAA^N$
be a closed immersion of codimension $r$. Let $F_1,\ldots,F_d$ be
general elements in the ideal of $I_X$ as in the above discussion.
Consider the subscheme $M$ of $\AAA^N$ defined by $F_1,\ldots,F_r$
and let $U_M$ be the open subset of ${\rm Cont}^e({\rm Jac}_X)_m$
that is contained in ${\rm Cont}^e({\rm Jac}_M)_m$. When we vary the
subsets with $r$ elements of $\{1,\ldots,d\}$, the corresponding
open subsets cover ${\rm Cont}^e({\rm Jac}_X)_m$. Therefore it is
enough to prove the two assertions in the proposition over $U_M$.

We claim that it is enough to prove i) and ii) for $M$. Indeed, if $\gamma\in U_M$ can be lifted to
$J_{m+e}(X)$, then in particular it can be lifted to $J_{m+e}(M)$.
If we know i) for $M$, it follows that $\gamma$ can be lifted to an
arc $\delta\in J_{\infty}(M)$. Lemma~\ref{lemma_reduction} implies
that $\delta$ lies in $J_{\infty}(X)$, hence we have i) for $X$.
Moreover, suppose that ii) holds for $M$, hence the projection
$$\beta\colon\psi_{m+1}^M({\rm Cont}^e({\rm Jac}_M))
\to \psi_m^M({\rm Cont}^e({\rm Jac}_M))$$ is piecewise trivial with
fiber $\AAA^n$. Again, Lemma~\ref{lemma_reduction} implies that the
restriction of $\beta$ over $U_M\cap\psi_m^M({\rm Cont}^e({\rm
Jac}_M))$ coincides with the restriction of $\alpha$ over
$U_M\cap\psi_m({\rm Cont}^e({\rm Jac}_X))$. Therefore $X$ also
satisfies ii).

We now prove  the proposition for a subscheme $M$ defined by a
regular sequence $F_1,\ldots,F_r$ ($M$ might not be reduced, but we do not
need this assumption anymore). Consider an element
$u=(u_1,\ldots,u_N)\in J_m(M)$, where all $u_i$ lie in
$k[t]/(t^{m+1})$ (for the matrix computations that will follow we
consider $u$ as a column vector). We denote by $\widetilde{u}_i\in
k\llbracket \rrbracket]]$ the lifting of $u_i$ that has degree $\leq m$. Our
assumption is that ${\rm ord}(F_i(\widetilde{u}))\geq m+1$ for every
$i$. An element in the fiber $(\psi_m^M)^{-1}(u)$ is an $N$--uple
$w=\widetilde{u}+t^{m+1}v$ where $v=(v_1,\ldots,v_N)\in (k\llbracket t\rrbracket)^N$,
such that $F_i(w)=0$ for every $i$. Using the Taylor expansion, we
get
\begin{equation}\label{eq_F}
F_i(w)=F_i(\widetilde{u})+t^{m+1}\cdot\sum_{j=1}^N\frac{\partial
F_i}{\partial x_j}(\widetilde{u})v_j+t^{2(m+1)}A_i(\widetilde{u},v),
\end{equation}
where each $A_i$ has all terms of degree $\geq 2$ in the $v_j$.
 We write $F$ and $A$ for the column vectors
$(F_1,\ldots,F_r)$ and $(A_1,\ldots,A_r)$, respectively.

Let $J(\widetilde{u})$ denote the Jacobian matrix $(\partial
F_i(\widetilde{u})/\partial x_j)_{i\leq r,j\leq N}$. Since $u$ lies
in ${\rm Cont}^e({\rm Jac}_M)_m$, all $r$--minors of this matrix
have order $\geq e$. Moreover, after taking a suitable open cover of
${\rm Cont}^e({\rm Jac}_M)_m$ and reordering the variables, we may
assume that the determinant of the submatrix $R(\widetilde{u})$ on
the first $r$ columns of $J(\widetilde{u})$ has order precisely $e$.
If $R^*(\widetilde{u})$ denotes the classical adjoint of the matrix
$R(\widetilde{u})$, then
$$R^*(\widetilde{u})\cdot J(\widetilde{u})=(t^e\cdot I_r, t^e\cdot J'(\widetilde{u}))$$
for some $r\times (N-r)$ matrix $J'(\widetilde{u})$. Indeed, for
every $i\leq r$ and $r+1\leq j\leq N$, the $(i,j)$ entry of
$R^*(\widetilde{u})\cdot J(\widetilde{u})$ is equal, up to a sign,
with the $r$--minor of $J(\widetilde{u})$ on the columns
$1,\ldots,i-1,i+1,\ldots,r,j$. Therefore its order is $\geq e$.

Since the determinant of $R^*(\widetilde{u})$ is nonzero, it follows
that $F(w)=0$ if and only if $R^*(\widetilde{u})\cdot F(w)=0$. By
equation (\ref{eq_F}) we have
\begin{equation}\label{eq_F2}
R^*(\widetilde{u})\cdot F(w)=R^*(\widetilde{u})\cdot
F(\widetilde{u})+t^{m+e+1}\cdot (I_r,J'(\widetilde{u}))\cdot v+
t^{2m+2}\cdot R^*(\widetilde{u})\cdot A(\widetilde{u},v).
\end{equation}
Note that since $m\geq e$ we have $2m+2>m+e+1$.

We claim that there is $v$ such that $F(\widetilde{u}+t^{m+1}v)=0$
if and only if
\begin{equation}\label{eq_F3}
\ord(R^*(\widetilde{u})\cdot F(\widetilde{u}))\geq m+e+1.
\end{equation}
Indeed, the fact that this condition is necessary follows
immediately from (\ref{eq_F2}). To see that it is also sufficient,
suppose that (\ref{eq_F3}) holds, and let us show that we can find
$v$ such that $F(\widetilde{u}+t^{m+1}v)=0$. We write
$v_i=\sum_{j}v_i^{(j)}t^j$ and determine inductively the
$v_i^{(j)}$. If we consider the term of order $m+e+1$ on the
right-hand side of (\ref{eq_F2}), then we see that we can choose
$v_{r+1}^{(0)},\ldots,v_N^{(0)}$ arbitrarily, and then the other
$v_i^{(0)}$ are uniquely determined. In the term of order
$t^{m+e+2}$, the contribution of the part coming from
$R^*(\widetilde{u})\cdot A(\widetilde{u},v)$ involves only the
$v_i^{(0)}$. It follows that again we may choose
$v_{r+1}^{(1)},\ldots,v_N^{(1)}$ arbitrarily, and then
$v^{(1)}_1,\ldots,v^{(1)}_r$ are determined uniquely such that the
coefficient of $t^{m+e+2}$ in $R^*(\widetilde{u})\cdot
F(\widetilde{u}+t^{m+1}v)$ is zero. Continuing this way we see that
we can find $v$ such that $F(\widetilde{u}+t^{m+1}v)=0$. This
concludes the proof of our claim.  Since the fiber over $u$ in
$\psi_{m+1}(J_{\infty}(M))$ corresponds to those
 $(v_1^{(0)},\ldots,v_N^{(0)})$ such that there is $v$ with
 $F(\widetilde{u}+t^{m+1}v)=0$, it follows from our description that
  this is a linear
 subspace of codimension $r$ of $\AAA^N$.

Note that if there is $v$ such that
$\ord\,F(\widetilde{u}+t^{m+1}v)\geq m+e+1$, then as above we get
that $\ord(R^*(\widetilde{u})\cdot F(\widetilde{u}))\geq m+e+1$. We
deduce that if $u$ can be lifted to $J_{m+e}(M)$, then $u$ can be
lifted to $J_{\infty}(M)$, which proves i). Moreover, the above
computation shows that over the set $W$ defined by (\ref{eq_F3}) in
our locally closed subset of $J_m(M)$, the inclusion
$$\psi_{m+1}({\rm Cont}^e({\rm Jac}_M))\subseteq
J_m(M)\times\AAA^N$$ is, at least set-theoretically, an affine
bundle with fiber $\AAA^{N-r}$. This proves ii) and completes the
proof of the proposition.
\end{proof}

\begin{remark}\label{remark_fiber2}
It follows from the above proof that the assertions of the
proposition hold also for a locally complete intersection scheme
(the scheme does not have to be reduced).
\end{remark}

\bigskip

We now discuss the fibers of the truncation maps between jet spaces
without restricting to the jets that can be lifted to the space of
arcs.

\begin{proposition}\label{fiber6}{\rm (}\cite{Lo}{\rm )}
Let $X$ be a scheme of finite type over $k$. For every nonnegative
integers $m$ and $p$, with $p\leq m\leq 2p+1$, consider the projection
$$\pi_{m,p}\colon J_m(X)\to J_p(X).$$
\begin{enumerate}
\item[i)] If $\gamma\in J_p(X)$ is such that
$\pi_{m,p}^{-1}(\gamma)$ is non-empty, then scheme-theoretically we
have
\begin{equation}\label{equation_fiber6}
\pi_{m,p}^{-1}(\gamma)\simeq {\rm
Hom}_{k[t]/(t^{p+1})}(\gamma^*\Omega_X,(t^{p+1})/(t^{m+1})).
\end{equation}
\item[ii)] Suppose that $X$ has pure dimension $n$
and that for $e=\ord_{\gamma}({\rm Jac}_X)$ we have $2p\geq m\geq
e+p$. If $X$ is either locally complete intersection or reduced,
and if $\pi_{m,p}^{-1}(\gamma)\neq\emptyset$,
then
$$\pi_{m,p}^{-1}(\gamma)\simeq\AAA^{e+(m-p)n}.$$
\end{enumerate}
\end{proposition}

\begin{proof}
Note that $\gamma$ corresponds to a ring homomorphism $\OO_{X,x}\to
k[t]/(t^{m+1})$, for some $x\in X$. Our assumption on $m$ and $p$
implies that $(t^{p+1})/(t^{m+1})$ is a $k[t]/(t^{p+1})$--module.
Therefore the right-hand side of (\ref{equation_fiber6}) is
well-defined. It is a finite-dimensional $k$--vector space, hence it
is an affine space.

In order to describe it, we use the structure of finitely generated
modules over $k\llbracket t\rrbracket$ to write a free presentation
$$(k[t]/(t^{p+1})^{\oplus N}\overset{A}\to
(k[t]/(t^{p+1}))^{\oplus N}\to \gamma^*\Omega_X\to 0,$$ where $A$ is
the diagonal matrix ${\rm diag}(t^{a_1},\ldots,t^{a_{N}})$, with
$0\leq a_1\leq\ldots\leq a_{N}\leq p+1$. In this case the right-hand
side of (\ref{equation_fiber6}) is isomorphic to $\AAA^{\ell}$,
where $\ell=\sum_i\min\{a_i,m-p\}$. Note also that its $R$--valued
points are in natural bijection with
$${\rm Der}_k(\OO_{X,x},t^{p+1}R[t]/t^{m+1}R[t]).$$

We first show  that it is enough to prove i). Suppose that we are in
the setting of ii).  We use the above description of the right-hand
side of (\ref{equation_fiber6}). It follows from the definition of
$e$ that $\sum_{i=1}^{N-n}a_i=e$. In particular, $a_i\leq e\leq
m-p$ for $i\leq N-n$. In order to deduce ii) from i) it is enough to show that
$a_i=p+1$ for $i>N-n$. If $\delta$ is an element in
$\pi_{m,p}^{-1}(\gamma)$, then by taking a free presentation of
$\delta^*\Omega_X$, we see that $A$ is the reduction mod $(t^{p+1})$
of a matrix ${\rm diag}(t^{b_1},\ldots,t^{b_N})$ with $0\leq
b_1\leq\ldots\leq b_N\leq m+1$. We have $a_i=b_i$ if $b_i\leq p$ and
$a_i=p+1$ otherwise. Either of our two conditions on $X$ implies
that ${\rm Fitt}^{n-1}(\Omega_X)=0$, hence
$${\rm ord}_{\delta}({\rm Fitt}^{n-1}(\Omega_X))
\geq m+1,$$ and therefore $b_1+\ldots+b_{N-n+1}\geq m+1$. We deduce
that for every $i\geq N-n+1$ we have $b_i\geq m+1-e\geq p+1$, hence
$a_i=p+1$.

Therefore it is enough to prove i). We may clearly assume that
$X={\rm Spec}(S)$ is affine. We start with the following
observation. If $\beta$ is an $R$--valued point in $J_p(X)$, then
either the fiber $\pi_{m,p}^{-1}(\beta)$ is empty, or it is a
principal homogeneous space over ${\rm
Der}_k(S,t^{p+1}R[t]/t^{m+1}R[t])$, where $t^{p+1}R[t]/t^{m+1}R[t]$
becomes an $S$--module via $\beta\colon S\to R[t]/(t^{p+1})$.
Indeed, if $D\in {\rm Der}_k(S,t^{p+1}R[t]/t^{m+1}R[t])$ and if
$\alpha\colon S\to R[t]/(t^{m+1})$ corresponds to an $R$--valued
point in $J_m(X)$ lying over $\beta$, then $\alpha+D$ gives another
$R$--valued point over $\beta$. Moreover, every other element in
$\pi_{m,p}^{-1}(\beta)$ arises in this way for a unique derivation
$D$.

We see that if $\delta$ is a fixed $k$--valued point in
$\pi_{m,p}^{-1}(\gamma)$, then we get a morphism
$${\rm
Hom}_{k[t]/(t^{p+1})}(\Omega_S\otimes_Sk[t]/(t^{p+1}),(t^{p+1})/(t^{m+1}))\to
\pi_{m,p}^{-1}(\gamma).$$ This is an isomorphism since it induces a
bijection at the level of $R$--valued points for every $R$.
\end{proof}

\begin{remark}\label{same_fiber6}
Let $X$ be a reduced scheme of pure dimension $n$. Suppose that $m$,
$p$ and $e$ are nonnegative integers such that $2p\geq m\geq e+p$
and $\gamma\in J_p(X)$ is such that $\ord_{\gamma}({\rm Jac}_X)=e$.
Assume also that $X$ is a closed subscheme of a locally complete
intersection scheme $M$ of the same dimension such that ${\rm
ord}_{\gamma}({\rm Jac}_M)=e$ (if $X$ is embedded in some $\AAA^N$,
then one can take $M$ to be generated by $(N-n)$ general elements in
the ideal of $X$). Consider the commutative diagram
\[
\begin{CD}
J_{m}(X)@>>> J_{m}(M)\\
@VV{\pi_{m,p}^X}V@VV{\pi_{m,p}^M}V\\
J_p(X) @>>>J_p(M)
\end{CD}
\]
where the horizontal maps are inclusions. It follows from
Proposition~\ref{fiber6} that the scheme-theoretic fibers of
$\pi_{m,p}^X$ and $\pi_{m,p}^M$ over $\gamma$ are equal.

Indeed, note first that if
$(\pi_{m,p}^M)^{-1}(\gamma)\neq\emptyset$, then $\gamma$ can be
lifted to $J_{\infty}(M)$ by Proposition~\ref{fiber2} (see also
Remark~\ref{remark_fiber2}). On the other hand such a lifting would
lie in $J_{\infty}(X)$ by Lemma~\ref{lemma_reduction}, hence
$(\pi_{m,p}^X)^{-1}(\gamma)\neq\emptyset$. In this case, it follows
from Proposition~\ref{fiber6} that both fibers are affine spaces of
the same dimension, one contained in the other, hence they are equal.
\end{remark}

\begin{remark}\label{remark_nonsingular_case}
Suppose that $X$ is a nonsingular variety of dimension $n$, and
suppose that $m\leq 2p+1$. On $J_p(X)$ we have a geometric vector
bundle $E$ whose fiber over $\gamma$ is ${\rm
Hom}_{k[t]/(t^{p+1})}(\gamma^*\Omega_X,(t^{p+1})/(t^{m+1}))$. If we
consider this as a group scheme over $J_p(X)$, then the argument in
the proof of Proposition~\ref{fiber6} shows that we have an action
of $E$ on $J_m(X)$ over $J_p(X)$. Moreover, whenever we have a
section of the projection of $\pi_{m,p}$ we get an isomorphism of
$J_m(X)$ with $E$. We can always find such a section if we restrict
to an affine open subset of $X$ on which $\Omega_X$ is trivial.
\end{remark}

We will need later the following global version of the assertion in
Proposition~\ref{fiber6} ii).

\begin{proposition}\label{fiber1}
Let $X$ be a scheme of pure dimension $n$ that is either reduced or
a locally complete intersection. If $m$, $p$ and $e$ are nonnegative
integers such that $2p\geq m\geq p+e$, then the canonical projection
$\pi_{m,p}\colon J_{m}(X)\to J_p(X)$ induces a piecewise trivial
fibration
$${\rm Cont}^e({\rm Jac}_X)_{m}\to {\rm Cont}^e({\rm Jac}_X)_p\cap
{\rm Im}(\pi_{m,p})$$ with fiber $\AAA^{(m-p)n+e}$.
\end{proposition}

\begin{proof}
We need to "globalize" the argument in the proof of
Proposition~\ref{fiber6}. Note that we may assume that $X$ is
locally a complete intersection. Indeed, we may assume first that
$X$ is affine. If $X$ is reduced, arguing as in the proof of
Proposition~\ref{fiber2} we may cover ${\rm Cont}^e({\rm Jac}_X)_p$
by open subsets $U_i$ such that there are $n$--dimensional locally
complete intersection schemes $M_i$ containing $X$, with
$$U_i\subseteq {\rm Cont}^e({\rm Jac}_{M_i})_p\subseteq J_p(M_i).$$
It follows from Remark~\ref{same_fiber6} that knowing the assertion
in the proposition for each $M_i$, we get it also for $X$.

Therefore we may assume that $X$ is a closed subscheme of $\AAA^N$
of codimension $r$, defined by $f_1,\ldots,f_r$. Write
$f=(f_1,\ldots,f_r)$, which we consider as a vertical vector.
Suppose that
$$u=(u_1,\ldots,u_N)\in {\rm Cont}^e({\rm Jac}_X)_p,$$ where
$u_i\in k[t]/(t^{p+1})$ for every $i$. We denote by
$\widetilde{u}\in (k[t]/(t^{m+1}))^N$ the lifting of $u$ having each
entry of degree $\leq p$. The fiber of $\pi_{m,p}$ over $u$ consists
of those $\widetilde{u}+t^{p+1}v$ such that
$f(\widetilde{u}+t^{p+1}v)=0$ in $(k[t]/(t^{m+1}))^N$. Here
$v=(v_1,\ldots,v_N)$ where $v_i=\sum_{j=0}^{m-p-1}v^{(j)}_it^j$.

Denote by $J(\widetilde{u})$ the Jacobian matrix $(\partial
f_i(\widetilde{u})/\partial x_j)_{i\leq r,j\leq N}$.
 Using the Taylor
expansion we see that
\begin{equation}\label{eq_F6}
f(\widetilde{u}+t^{p+1}v)=f(\widetilde{u})+t^{p+1}\cdot
J(\widetilde{u})v
\end{equation}
(there are no further terms since $2(p+1)\geq m+1$).

Note that by assumption we can write
$f(\widetilde{u})=t^{p+1}g(u)$ where
$g(u)=\,\,\,\,\,\,\,\,\,\,\,\,\,\,\,\,$
$(\sum_{j=0}^{m-p-1}g_{i,j}(u)t^j)_i$. If we denote by
$\overline{J}(u)$ the reduction of $J(\widetilde{u})$ mod $t^{m-p}$,
we see that the condition on $v$ becomes
\begin{equation}\label{eq_F7}
-g(u)=\overline{J}(u)\cdot v,
\end{equation}
where the equality is in $(k[t]/(t^{m-p}))^{r}$.

It follows from the structure theory of matrices over principal
ideal domains, applied to a lifting of $\overline{J}(u)$ to a matrix
over $k\llbracket t\rrbracket$, that we can find invertible matrices $A$ and $B$ over
$k[t]/(t^{m-p})$ such that $A\cdot \overline{J}(u)\cdot B=\left({\rm
diag}(t^{a_1},\ldots,t^{a_r}),{\mathbf {0}}\right)$, with $0\leq a_i\leq m-p$. 
Moreover, after partitioning ${\rm
Cont}^e({\rm Jac}_X)_p$ into suitable locally closed subsets, we may
assume that the $a_i$ are independent of $u$ and that $A=A(u)$ and
$B=B(u)$, where the entries of $A(u)$ and $B(u)$ are regular
functions of $u$.

 Since the ideal generated by the $r$--minors of $\overline{J}(u)$
is $(t^e)$, we see that $a_1+\ldots+a_r=e$. If we write $A(u)\cdot
g(u)=(h_1(u),\ldots,h_r(u))$, we see that $u$ lies in
the image of $\pi_{m,p}$ if and only if $\ord(h_i(u))\geq a_i$ for
every $i\leq r$. Moreover, if we put $v'=B(u)^{-1}v$, then we see
that our condition gives the values of $t^{a_i}v'_i$ for $i\leq r$.
Therefore the set of possible $v$ is isomorphic to an affine space
of dimension $(N-r)(m-p)+\sum_{i=1}^ra_i=n(m-p)+e$. Since the
equations defining the fiber over $u$ depend algebraically on $u$,
we get the assertion of the proposition.
\end{proof}

\section{Cylinders in spaces of arcs}

We start by giving some applications of Proposition~\ref{fiber2}.
For every scheme $X$, a \emph{cylinder} in $J_{\infty}(X)$ is a
subset of the form $C=\psi_m^{-1}(S)$, for some $m$ and some
constructible subset $S\subseteq J_m(X)$. From now on, unless
explicitly mentioned otherwise, we assume  that $X$ is reduced and
of pure dimension $n$.

\begin{lemma}\label{cylinder1}
If $C\subseteq J_{\infty}(X)$ is a cylinder, then $C$ is not thin if
and only if it is not contained in $J_{\infty}(X_{\rm sing})$.
\end{lemma}

\begin{proof}
We need to show that for every closed subset $Y$ of $X$ with
$\dim(Y)<\dim(X)$, and every cylinder $C\not\subseteq
J_{\infty}(X_{\rm sing})$, we have $C\not\subseteq J_{\infty}(Y)$.
If this is not the case, then arguing by Noetherian induction we may
choose a minimal $Y$ for which there is a cylinder $C\not\subseteq
J_{\infty}(X_{\rm sing})$ with $C\subseteq J_{\infty}(Y)$. After
replacing $C$ by a suitable $C\cap {\rm Cont}^e({\rm Jac}_X)$, we
may assume that $C\subseteq {\rm Cont}^e({\rm Jac}_X)$. It follows
from Proposition~\ref{fiber2} that if $m\gg 0$, then the maps
$\psi_{m+1}(C)\to \psi_m(C)$ are piecewise trivial, with fiber
$\AAA^n$.

Note that $Y$ has to be irreducible. Indeed, if $Y=Y_1\cup Y_2$,
with $Y_1$ and $Y_2$ both closed and different from $Y$, then either
$C\cap J_{\infty}(Y_1)$ or $C\cap J_{\infty}(Y_2)$ is not contained
in $J_{\infty}(X_{\rm sing})$. This contradicts the minimality of
$Y$.

Using again the fact that $Y$ is minimal, we see that
$C\not\subseteq J_{\infty}(Y_{\rm sing})$ (we consider $Y$ with the
reduced structure). After replacing $C$ with some $C\cap {\rm
Cont}^{e'}({\rm Jac}_Y)$, we may assume that $C\subseteq {\rm
Cont}^{e'}({\rm Jac}_Y)$. Since $C$ is a cylinder also in
$J_{\infty}(Y)$, it follows from Proposition~\ref{fiber2} that if
$m\gg 0$, then the projection $\psi_{m+1}(C)\to \psi_m(C)$ is
piecewise trivial with fiber $\AAA^{\dim(Y)}$. This is a
contradiction, and completes the proof of the lemma.
\end{proof}

\begin{corollary}\label{cor_cylinder1}
Let $f\colon X'\to X$ be a proper birational morphism of reduced,
pure-dimensional schemes. If $\gamma\in
\psi_m^X(J_{\infty}(X)\smallsetminus J_{\infty}(X_{\rm sing}))$,
then $\gamma$ lies in the image of $f_m$.
\end{corollary}

\begin{proof}
If $C=(\psi^X_m)^{-1}(\gamma)$, then $C$ is a cylinder that is not
contained in $J_{\infty}(X_{\rm sing})$. Let $Z\subset X$ be a
closed subset with $\dim(Z)<\dim(X)$ such that $f$ is an isomorphism
over $X\smallsetminus Z$. It follows from Proposition~\ref{prop11}
that $J_{\infty}(X)\smallsetminus J_{\infty}(Z)\subseteq {\rm
Im}(f_{\infty})$. Since $C\not\subseteq J_{\infty}(Z)$ by the lemma,
we deduce that there is $\delta\in J_{\infty}(X')$ such that
$f_m(\psi^{X'}_m(\delta))=\psi_m^X(f_{\infty}(\delta))=\gamma$.
\end{proof}

\begin{corollary}\label{cor2_cylinder1}
Let $f$ be as in the previous corollary, with $X'$ nonsingular. If
$k$ is uncountable, then $J_{\infty}(X)\smallsetminus
J_{\infty}(X_{\rm sing})\subseteq {\rm Im}(f_{\infty})$.
\end{corollary}

\begin{proof}
Let $\gamma\in J_{\infty}(X)\smallsetminus J_{\infty}(X_{\rm
sing})$. It follows from Corollary~\ref{cor_cylinder1} that for
every $m$ we have $\gamma_m:=\psi_m^X(\gamma)\in {\rm Im}(f_m)$.
Therefore we get a decreasing sequence
$$\cdots\supseteq(\psi_m^{X'})^{-1}(f_m^{-1}(\gamma_m))\supseteq
(\psi_{m+1}^{X'})^{-1}(f_{m+1}^{-1}(\gamma_{m+1}))\supseteq\cdots$$
of nonempty cylinders. Lemma~\ref{sequence} below implies that there
is $\delta$ in the intersection of all these cylinders. Therefore
$\psi_m^X(f_{\infty}(\delta))=\gamma_m$ for all $m$, hence
$\gamma=f_{\infty}(\delta)$.
\end{proof}

\begin{lemma}\label{sequence}{\rm (}\cite{batyrev}{\rm )}
If $X$ is nonsingular and $k$ is uncountable, then every decreasing
sequence of cylinders
$$C_1\supseteq\cdots \supseteq C_m\supseteq\cdots$$
has nonempty intersection.
\end{lemma}

\begin{proof}
Since the projections $\psi_m$ are surjective, it follows from
Chevalley's Constructibility Theorem that the image of every
cylinder in $J_m(X)$ is constructible. Consider the decreasing
sequence
$$\psi_0(C_1)\supseteq \psi_0(C_2)\supseteq\cdots$$
of nonempty constructible subsets. Since $k$ is uncountable, the
intersection of this sequence is nonempty. Let $\gamma_0$ be an
element in this intersection.

Since $\gamma_0$ lies in the image of every $C_m$, we see that all
the constructible subsets  in the decreasing sequence 
$$\psi_1(C_1)\cap\pi_{1,0}^{-1}(\gamma_0)\supseteq
\psi_1(C_2)\cap\pi_{1,0}^{-1}(\gamma_0)\supseteq\cdots$$ are
nonempty. Therefore there is $\gamma_1$ contained in their
intersection. Continuing in this way we get $\gamma_m\in J_m(X)$ for
every $m$ such that $\pi_{m,m-1}(\gamma_m)=\gamma_{m-1}$ for every
$m$ and $\gamma_m\in\psi_m(C_p)$ for every $p$. Therefore
$(\gamma_m)_m$ determines an arc $\gamma\in J_{\infty}(X)$ whose image
in  $J_m(X)$ is equal to $\gamma_m$. Since each $C_p$ is a cylinder and
$\psi_m(\gamma) \in\psi_m(C_p)$ for every $m$, we see that
$\gamma\in C_p$. Hence $\gamma\in\cap_{p\geq 1}C_p$.
\end{proof}

\begin{remark}\label{remark_sequence}
Note that in the above lemma, the hypothesis that $X$ is nonsingular
was used only to ensure that the image in $J_m(X)$ of a cylinder is
a constructible set. We will prove this below for an arbitrary
scheme $X$ (see Corollary~\ref{image_cylinder}), and therefore the
lemma will hold in this generality.
\end{remark}

\begin{remark}
If ${\rm char}(k)=0$, then the assumption that $X'$ is nonsingular
is not necessary in Corollary~\ref{cor2_cylinder1}. Indeed, we can
take a resolution of singularities $g\colon X''\to X'$ and we
clearly have ${\rm Im}(f\circ g)_{\infty}\subseteq {\rm
Im}(f_{\infty})$.
\end{remark}

\smallskip

We have seen in Proposition~\ref{fiber2} that for a reduced
pure-dimensional scheme $X$ the set $\psi_m({\rm Cont}^e({\rm
Jac}_X))$ is constructible. In fact, the image of every cylinder is
constructible, as follows from the following  result of
Greenberg.

\begin{proposition}\label{image_constructible}{\rm (}\cite{greenberg}{\rm )}
 For an arbitrary scheme $X$ and every $m$, the image of
$J_{\infty}(X)\to J_m(X)$ is constructible.
\end{proposition}

\begin{proof}
We give the proof assuming that ${\rm char}(k)=0$. For a proof in
the general case, see \cite{greenberg}. We do induction on
$\dim(X)$, the case $\dim(X)=0$ being trivial. If $X_1,\ldots,X_r$
are the irreducible components of $X$, with the reduced structure,
then $J_{\infty}(X)=J_{\infty}(X_1)\cup\ldots \cup J_{\infty}(X_r)$.
Hence the image of $J_{\infty}(X)$ is equal to the union of the
images of the $J_{\infty}(X_i)$ in $J_m(X_i)\subseteq J_m(X)$.
Therefore we may assume that $X$ is reduced and irreducible.

Let $f\colon X'\to X$ be a resolution of singularities. Since $X'$
is nonsingular, the projection $J_{\infty}(X')\to J_m(X')$ is
surjective, hence ${\rm Im}(f_m)\subseteq {\rm Im}(\psi_m^X)$.
Moreover,  Corollary~\ref{cor_cylinder1} gives
$\psi_m^X(J_{\infty}(X) \smallsetminus J_{\infty}(X_{\rm
sing}))\subseteq {\rm Im}(f_m)$. Therefore
$$\psi_m^X(J_{\infty}(X))={\rm Im}(f_m)\cup \psi_m^X(J_{\infty}(X_{\rm
sing})).$$ The first term on the right-hand side is constructible by
Chevalley's Constructibility Theorem, while the second term is
constructible by induction. This implies that
$\psi_m^X(J_{\infty}(X))$ is constructible.
\end{proof}

\begin{corollary}\label{image_cylinder}
For an arbitrary scheme $X$, the image of a cylinder $C$ by the
projection $J_{\infty}(X)\to J_m(X)$ is constructible.
\end{corollary}

\begin{proof} Let
$C=\psi_p^{-1}(A)$, where $A\subseteq J_p(X)$ is constructible. If
$m\geq p$, then $\psi_m(C)=\psi_m(J_{\infty}(X))\cap
\pi_{m,p}^{-1}(A)$, hence it is constructible by the proposition.
The constructibility for $m<p$ now follows  from Chevalley's Theorem.
\end{proof}

Proposition~\ref{image_constructible} is deduced in \cite{greenberg}
from the fact that for every $m$ there is $p\geq m$ such that the
image of the projection $\psi_m\colon J_{\infty}(X)\to J_m(X)$ is
equal to the image of $\pi_{p,m}\colon J_p(X)\to J_m(X)$ (in fact,
Greenberg also shows that one can take $p=L(m)$ for a suitable
linear function $L$). We now  show that if we assume $k$ uncountable,
then this follows from the above proposition.

\begin{corollary}\label{cor_image_constructible}
If $k$ is uncountable, then for an arbitrary scheme $X$ and every
$m$ there is $p\geq m$ such that the image of $\psi_m$ is equal to
the image of $\pi_{p,m}$.
\end{corollary}

\begin{proof}
Since ${\rm Im}(\psi_m)$ is constructible by
Proposition~\ref{image_constructible} and each ${\rm Im}(\pi_{p,m})$
is constructible by Chevalley's Theorem, the assertion follows if we
show
\begin{equation}\label{eq_intersection}
{\rm Im}(\psi_m)=\bigcap_{p\geq m}{\rm Im}(\pi_{p,m})
\end{equation}
(we use the fact that $k$ is uncountable). The inclusion
"$\subseteq$" is obvious. For the reverse inclusion we argue as in
the proof of Lemma~\ref{sequence} to show that if
 $\gamma_m\in\cap_{p\geq m} {\rm Im}(\pi_{p,m})$,
then we can find $\gamma_q\in J_q(X)$ for every $q\geq m+1$ such
that $\pi_{q,q-1}(\gamma_{q})=\gamma_{q-1}$. The sequence
$(\gamma_q)_q$ defines an element $\gamma\in J_{\infty}(X)$ lying
over $\gamma_m$.
\end{proof}

We give one more result about the fibers of the truncation maps
between the images of the spaces of arcs (one should compare this
with Proposition~\ref{fiber2}).

\begin{proposition}\label{fiber3}{\rm (}\cite{DL}{\rm )}
If $X$ is a scheme of dimension $n$, then for every $m\geq p$, all
fibers of the truncation map
$$\phi_{m,p}\colon\psi_{m}(J_{\infty}(X))\to\psi_p(J_{\infty}(X))$$
have dimension $\leq (m-p)n$.
\end{proposition}

\begin{proof}
Note that the sets in the statement are constructible by
Proposition~\ref{image_constructible}. Clearly, it is enough to
prove the proposition when $m=p+1$. We may assume that $X$ is a
closed subscheme of $\AAA^N$ defined by $F_1,\ldots,F_r$. Consider
$\gamma_p\in J_p(X)$ given by $u=(u_1,\ldots,u_N)$ where $u_i\in
k[t]$ with $\deg(u_i)\leq p$.

Let $T={\rm Spec}\,k[t]$. Consider the subscheme ${\mathcal Z}$ of
$T\times\AAA^N$ defined by $I_{\mathcal Z}=(F_1(u+t^{p+1}x),\ldots,
F_r(u+t^{p+1}x))$. We have a subscheme ${\mathcal Z}'\subseteq
{\mathcal Z}$ defined by
$$I_{{\mathcal Z}'}=(f\mid hf\in I_{\mathcal Z}\,{\rm for}\,{\rm
some}\,{\rm nonzero}\,h\in k[t]).$$ Note that by construction
${\mathcal Z}'$ is flat over $T$, and ${\mathcal Z}={\mathcal Z}'$
over the generic point of $T$.

The generic fiber of ${\mathcal Z}$ over $T$ is isomorphic to
$X\times_kk(t)$. Since ${\mathcal Z}'$ is flat over $T$, it follows
that the fiber of ${\mathcal Z}'$ over the origin is either empty or
has dimension $n$.

On the other hand, an element in the fiber of $\phi_{p+1,p}$ over
$\gamma_p$ is the $(p+1)$--jet of an arc in $X$ given by
$u+t^{p+1}w$ for some $w\in (k\llbracket t\rrbracket)^N$. Since $F_i(u+t^{p+1}w)=0$
for every $i$, it follows from the definition of $I_{{\mathcal Z}'}$
that if
 $f\in I_{{\mathcal Z}'}$, then $f(t,w)=0$. Hence the fiber
 of $\phi_{p+1,p}$ over $\gamma_p$ can be embedded in
the fiber of ${\mathcal Z}'$ over the origin, and its dimension is
$\leq n$.
\end{proof}

\smallskip

We now 
discuss  the notion of codimension for cylinders in spaces of
arcs. In the remaining part of this section we assume that $k$ is
uncountable, and also that ${\rm char}(k)=0$ (this last condition is
due only to the fact that we use resolutions of singularities).

Let $X$ be a scheme of pure dimension $n$ that is either reduced or
locally complete intersection, and let $C=\psi_p^{-1}(A)$ be a
cylinder, where $A$ is a constructible subset of $J_p(X)$. If
$C\subseteq {\rm Cont}^e({\rm Jac}_X)$ and $m\geq \max\{p,e\}$, then
we put $\codim(C):=(m+1)n-\dim(\psi_m(C))$. We refer to \S 9.1
for a quick review of some basic facts about the dimension of
constructible subsets. Note that by Proposition~\ref{fiber2} (see
also Remark~\ref{remark_fiber2}), this is well-defined. Moreover, it
is a nonnegative integer:  by Theorem~\ref{thm1}, the closure of
$\psi_m(J_{\infty}(X))$ is equal to the closure in $J_m(X)$ of the
$m^{\rm th}$ jet scheme of the nonsingular locus of $X_{\rm red}$.
Therefore it is a set of pure dimension $(m+1)n$ (the fact that
$\dim\,\psi_m(J_{\infty}(X))=(m+1)n$ follows also from
Proposition~\ref{fiber3}).

For an arbitrary cylinder $C$ we put $C^{(e)}:=C\cap {\rm
Cont}^e({\rm Jac}_X)$ and
$$\codim(C):=\min\{\codim(C^{(e)})\mid e\in\NN\}$$
(by convention, if $C\subseteq J_{\infty}(X_{\rm sing})$, we have
$\codim(C)=\infty$). It is clear that if $C_1$ and $C_2$ are
cylinders, then $\codim(C_1\cup
C_2)=\min\{\codim(C_1),\codim(C_2)\}$. In particular, if
$C_1\subseteq C_2$, then $\codim(C_1)\geq\codim(C_2)$.

\begin{proposition}\label{countable_union}
Suppose that $X$ is reduced and let $C$ be a cylinder in
$J_{\infty}(X)$. If we have disjoint cylinders $C_i\subseteq C$ for
$i\in\NN$ such that the complement $C\smallsetminus
\bigsqcup_{i\in\NN}C_i$ is thin, then
$\lim_{i\to\infty}\codim(C_i)=\infty$ and
$\codim(C)=\min_i\codim(C_i)$.
\end{proposition}

\noindent Note that the proposition implies that for every cylinder $C$ we
have
$$\lim_{e\to\infty}\codim(C^{(e)})=\infty.$$
We will prove Proposition~\ref{countable_union} at the same time
with the following proposition.

\begin{proposition}\label{limit}
If $X$ is reduced and $Y$ is a closed subscheme of $X$ with
$\dim(Y)<\dim(X)$, then
$$\lim_{m\to\infty}\codim({\rm Cont}^{\geq m}(Y))=\infty.$$
\end{proposition}

We first show  that these results hold when $X$ is nonsingular.
Let us start by making some comments about this special case.
Suppose for the moment that $X$ is nonsingular of pure dimension
$n$. Since the projections $J_{m+1}(X)\to J_m(X)$ are locally
trivial with fiber $\AAA^n$, cylinders are much easier to understand
in this case. We say that a cylinder $C=\psi_m^{-1}(S)$ is
\emph{closed}, \emph{locally closed} or \emph{irreducible} if $S$ is
(the definition does not depend on $m$ by the local triviality of
the projection). Moreover, if $S$ is closed and $S=S_1\cup\ldots\cup
S_r$ is the irreducible decomposition of $S$, then we get a unique
decomposition into maximal irreducible closed cylinders
$C=\psi_m^{-1}(S_1)\cup\ldots\cup\psi_m^{-1}(S_r)$. The cylinders
$\psi_m^{-1}(S_i)$ are the \emph{irreducible components} of $C$.

Note that if $C=\psi_m^{-1}(S)$, then by definition
$\codim(C)=\codim(S,J_m(X))$. If $C\subseteq C'$ are closed
cylinders with $\codim(C)=\codim(C')$, then every irreducible
component of $C$ whose codimension is equal to $\codim(C)$ is also
an irreducible component of $C'$.

\begin{proof}[Proof of Propositions~\ref{countable_union} and \ref{limit}]
We start by noting that if Proposition~\ref{limit} holds on $X$,
then Proposition~\ref{countable_union} holds on $X$, too. Indeed,
suppose that $\bigsqcup_{i\in\NN}C_i\subseteq C$, where all $C_i$ and
$C$ are cylinders, and that $C\smallsetminus\bigsqcup_iC_i$ is
contained in $J_{\infty}(Y)$, where $\dim(Y)<\dim(X)$. For every $m$
we have
$$C\subseteq {\rm Cont}^{\geq
m}(Y)\cup\bigcup_{i\in\NN}C_i.$$ It follows from
Lemma~\ref{sequence} that there is an integer $i(m)$ such that
\begin{equation}\label{eq_prop_cylinders}
C\subseteq{\rm Cont}^{\geq
m}(Y)\cup\bigcup_{i\leq i(m)}C_i.
\end{equation}
 In particular, for every
$i>i(m)$ we have $C_i\subseteq {\rm Cont}^{\geq m}(Y)$, hence
$\codim(C_i)\geq\codim\,{\rm Cont}^{\geq m}(Y)$. If
Proposition~\ref{limit} holds on $X$, it follows that
$$\lim_{i\to\infty}\codim(C_i)=\infty.$$

The second assertion in Proposition~\ref{countable_union} follows,
too. Indeed, note first that if all $C_i\subseteq J_{\infty}(X_{\rm
sing})$, then $C\subseteq J_{\infty}(Y\cup X_{\rm sing})$. Therefore
$C\subseteq J_{\infty}(X_{\rm sing})$ by Lemma~\ref{cylinder1}, and
the assertion is clear in this case. If $C_i\not\subseteq
J_{\infty}(X_{\rm sing})$ for some $i$, then $\codim(C)<\infty$. The
assertion in Proposition~\ref{limit} implies that there is $m$ such
that $\codim\,{\rm Cont}^{\geq m}(Y)>\codim(C)$. We deduce from
(\ref{eq_prop_cylinders}) that
$$\codim(C)\geq\min\{\codim(C_0),\ldots,\codim(C_{i(m)}),\codim\,{\rm
Cont}^{\geq m}(Y)\}.$$ Therefore $\codim(C)\geq\min_i\codim(C_i)$
and the reverse inequality is trivial.

We now prove Proposition~\ref{limit} when $X$ is nonsingular. We
have a decreasing sequence of closed cylinders $\{{\rm Cont}^{\geq
m}(Y)\}_{m\in\NN}$. Since
$$\codim\,{\rm Cont}^{\geq m}(Y)\leq\codim\,{\rm Cont}^{\geq
m+1}(Y)$$ for every $m$, it follows that if the limit in the
proposition is not infinity, then there is $m_0$ such that
$\codim\,{\rm Cont}^{\geq m}(Y)=\codim\,{\rm Cont}^{\geq m_0}(Y)$
for every $m\geq m_0$. Hence for all such $m$, the irreducible
components of ${\rm Cont}^{\geq m+1}(Y)$ of minimal codimension are
also components of ${\rm Cont}^{\geq m}(Y)$. It is easy to see that
this implies that there is an irreducible component $C$ of all ${\rm
Cont}^{\geq m}(Y)$ for $m\geq m_0$. Therefore $C\subseteq
J_{\infty}(Y)$, which contradicts Lemma~\ref{cylinder1}. By our
discussion at the beginning of the proof we see that both
propositions hold on nonsingular varieties.

In order to complete the proof it is enough to show that
Proposition~\ref{limit} holds for an arbitrary reduced
pure-dimensional scheme $X$. Let $f\colon X'\to X$ be a resolution
of singularities of $X$ (in other words $X'$ is the disjoint union
of resolutions of the irreducible components of $X$). Since
$f_{\infty}^{-1}({\rm Cont}^{\geq m}(Y))={\rm Cont}^{\geq
m}(f^{-1}(Y))$ and since we know that Proposition~\ref{limit} holds
on $X'$, we see that it is enough to prove that for every cylinder
$C\subseteq J_{\infty}(X)$, we have
$\codim(f_{\infty}^{-1}(C))\leq\codim(C)$.

We clearly have $\bigsqcup_{e\in\NN}f_{\infty}^{-1}(C^{(e)})\subseteq
f_{\infty}^{-1}(C)$ and the complement of this union is contained in
$J_{\infty}(f^{-1}(X_{\rm sing}))$. Since
Proposition~\ref{countable_union} holds on $X'$, we see that
$\codim(f_{\infty}^{-1}(C))=\min_e\codim\,f_{\infty}^{-1}(C^{(e)})$.
Therefore we may assume that $C=C^{(e)}$ for some $e$. In this case,
if $m\gg 0$, then
$$\codim(C)=(m+1)\dim(X)-\dim\,\psi^X_m(C)\geq
(m+1)\dim(X)-\dim\,f_m^{-1}(\psi^X_m(C))$$
$$=\codim\,(f_m\circ\psi_m^{X'})^{-1}(\psi_m^X(C))=\codim\,f_{\infty}^{-1}(C).$$
We have used the fact that $\psi^X_m(C)\subseteq {\rm Im}(f_m)$ by
Corollary~\ref{cor_cylinder1}. This completes the proof of the two
propositions.
\end{proof}

\begin{example}
Let $Z\subseteq\AAA^2$ be the curve defined by $x^2-y^3=0$. The
Jacobian ideal of $Z$ is ${\rm Jac}_Z=(x,y^2)$. Let $\pi\colon
J_{\infty}(Z)\to Z$ be the projection map. If $z\in Z$ is different
from the origin, then $z$ is a smooth
 point of $Z$ and $\codim(\pi^{-1}(z))=1$. On the other hand,
 if $z$ is the origin, then we can decompose $C=\pi^{-1}(z)$
 as
 $$C=J_{\infty}(z)\bigsqcup\left(\bigsqcup_{e>0}\{(u(t),v(t))\mid
 u(t)^2=v(t)^3, \ord\,u(t)=3e,\,\ord\,v(t)=2e\}\right).$$
 Note that the set corresponding to $e$ is precisely $C^{(3e)}$.
If we take $m=3e$, we see that $\psi_m(C^{(3e)})$ is equal to
$$\{(at^{3e}, b_0t^{2e}+\ldots+b_{e}t^{3e}\mid a^2=b_0^3, a\neq 0,
b_0\neq 0\}.$$ Therefore $\codim(C^{(3e)})=(3e+1)-(e+1)=2e$ for every $e\geq 1$, and
$\codim(C)=2$. Note that in this case the codimension of the special
fiber of $\pi$ is larger than that of the general fiber (compare
with the behavior of dimensions of fibers of morphisms of algebraic
varieties).
\end{example}

Proposition~\ref{countable_union} is a key ingredient in setting up
motivic integration (see \cite{batyrev} and \cite{DL}). We describe
one elementary application of this proposition to the definition of
another invariant of a cylinder, the "number of components of
minimal codimension".

Let $X$ be a reduced pure-dimensional scheme and $C$ a cylinder in
$J_{\infty}(X)$. If $C\subseteq {\rm Cont}^e({\rm Jac}_X)$, then we
take $m\gg 0$ and define $|C|$ to be the number of irreducible
components of $\psi_m(C)$ whose codimension is $\codim(C)$. Note
that by Proposition~\ref{fiber2} this number is independent of $m$.
For an arbitrary $C$, we put $|C|:=\sum_{e\in\NN}|C^{(e)}|$, where
the sum is over those $e$ such that $\codim(C^{(e)})=\codim(C)$
(Proposition~\ref{countable_union} implies that this is a finite
sum). With this definition, we see that under the hypothesis of
Proposition~\ref{countable_union} we have $|C|=\sum_i|C_i|$, the sum
being over the finite set of those $i$ with $\codim(C_i)=\codim(C)$.

If $X$ is a nonsingular variety and $C$ is a closed cylinder in
$J_{\infty}(X)$, then $|C|$ is equal to the number of irreducible
components of $C$ of minimal codimension.

\section{The Birational Transformation Theorem}

We now present  the fundamental result of the theory. Suppose that
$f\colon X'\to X$ is a proper birational morphism, with $X'$
nonsingular and $X$ reduced and of pure dimension $n$. The
Birational Transformation Theorem shows that in this case
$f_{\infty}$ induces at finite levels weakly piecewise trivial
fibrations.

The dimension of the fibers of these fibrations depends on the order
of vanishing along the \emph{Jacobian ideal} ${\rm Jac}_f$ of $f$. Consider the
morphism induced by pulling-back $n$--forms
$$f^*\Omega_X^n\to\Omega_{X'}^n.$$
Since $X'$ is nonsingular, $\Omega_{X'}^n$ is locally free of rank
one, hence the image of the above morphism can be written as ${\rm
Jac}_f\otimes\Omega_{X'}^n$ for a unique ideal ${\rm Jac}_f$ of
$\OO_{X'}$. In other words, we have ${\rm Jac}_f={\rm
Fitt}^0(\Omega_{X'/X})$.

If $X$ is nonsingular, too, then ${\rm Jac}_f$ is locally
principal, and it defines a subscheme supported on the exceptional
locus of $f$. In this case,  Proposition~\ref{prop11} implies that
 $f_{\infty}$ is injective on $J_{\infty}(X')\smallsetminus
J_{\infty}(V({\rm Jac}_f))$. In general, we have the following.

\begin{lemma}\label{injectivity2}
If $f\colon X'\to X$ is a proper birational morphism, with $X'$
nonsingular and $X$ reduced and pure-dimensional, and if $\gamma$,
$\gamma' \in J_{\infty}(X')$ are such that
$$\gamma\not\in J_{\infty}(V({\rm Jac}_f))\cup
f_{\infty}^{-1}(J_{\infty}(X_{\rm sing}))$$ and $f_{\infty}(\gamma)
=f_{\infty}(\gamma')$, then $\gamma=\gamma'$.
\end{lemma}

\begin{proof}
We argue as in the proof of Proposition~\ref{prop11}.  Since $f$ is
separated, it is enough to show that if $j\colon{\rm
Spec}\,k\llparenthesis t\rrparenthesis
\to {\rm Spec}\,k\llbracket t\rrbracket$ corresponds to $k\llbracket t\rrbracket\subset
k\llparenthesis t\rrparenthesis$, then $\gamma\circ j=\gamma'\circ j$.

Note that $U:=f^{-1}(X_{\rm reg})\smallsetminus V({\rm Jac}_f)$ is
an open subset of $X'$ that is the inverse image of an open subset
of $X$. Moreover, $f$ is invertible on $U$. By assumption,
$\gamma\circ j$ factors through $U$ and $f\circ\gamma\circ
j=f\circ\gamma'\circ j$. Therefore $\gamma'\circ j$ also factors
through $U$ and $\gamma\circ j=\gamma'\circ j$.
\end{proof}

\begin{theorem}\label{change_of_variable}
Let $f\colon X'\to X$ be a proper birational morphism, with $X'$
nonsingular and $X$ reduced and of pure dimension $n$. For
nonnegative integers $e$ and $e'$, we put
$$C_{e,e'}:={\rm Cont}^e({\rm Jac}_f)\cap f_{\infty}^{-1}({\rm
Cont}^{e'}({\rm Jac}_X)).$$ Fix $m\geq\max\{2e, e+e'\}$.
\begin{enumerate}
\item[i)] $\psi^{X'}_m(C_{e,e'})$ is a union of fibers of $f_m$.
\item[ii)] $f_m$ induces a weakly piecewise trivial fibration
with fiber $\AAA^e$
$$\psi^{X'}_m(C_{e,e'})\to f_m(\psi^{X'}_m(C_{e,e'})).$$
\end{enumerate}
\end{theorem}

In the case when also $X$ is nonsingular, this theorem is due to
Kontsevich \cite{kontsevich}. The case of singular $X$ is due to
Denef and Loeser \cite{DL}, while the proof we give below follows
\cite{Lo}. Note that in these references one makes the assumption
that the base field has characteristic zero, and therefore one gets
piecewise trivial fibrations in ii) above. For a version in the
context of formal schemes, allowing also positive characteristic,
but with additional assumptions on the morphism, see \cite{Se}. The
above theorem is at the heart of the Change of Variable Formula in
motivic integration (see \cite{batyrev}, \cite{DL}, and also \cite{Loeser}).

We start with some preliminary remarks. Let $f$ be as in the
theorem, and suppose that $\alpha\in J_{\infty}(X')$, with
$\ord_{\alpha}({\rm Jac}_f)=e$ and $\ord_{f_{\infty}(\alpha)}({\rm
Jac}_X)=e'$. Pulling-back via $\alpha$ the right exact sequence of
sheaves of differentials associated to $f$, we get an exact sequence
of $k\llbracket t\rrbracket$--modules
$$
\alpha^*(f^*\Omega_X)\overset{h}\to\alpha^*\Omega_{X'}\to\alpha^*\Omega_{X'/X}\to
0.
$$
 By assumption ${\rm Fitt}^0(\alpha^*\Omega_{X'/X})=(t^e)$, hence
$$\alpha^*(\Omega_{X'/X})\simeq k[t]/(t^{a_1})\oplus\ldots\oplus
k[t]/(t^{a_n})$$ for some $0\leq a_1\leq\ldots\leq a_n$ with
$\sum_ia_i=e$.

It follows that if $T={\rm Im}(h)$, then $T$ is free of rank $n$,
and in suitable bases of $T$ and $\alpha^*\Omega_{X'}$, the induced
map $g\colon T\to \alpha^*\Omega_{X'}$ is given by the diagonal
matrix with entries $t^{a_1},\ldots,t^{a_n}$. We get a decomposition
$\alpha^*(f^*\Omega_X)\simeq T\oplus {\rm Ker}(h)$, and therefore
$${\rm Fitt}^0({\rm Ker}(h))={\rm
Fitt}^n(\alpha^*(f^*\Omega_X))=(t^{e'}).$$ Hence ${\rm Ker}(h)\simeq
k[t]/(t^{b_1})\oplus\ldots\oplus k[t]/(t^{b_r})$ for some $0\leq
b_1\leq\ldots\leq b_r$ with $\sum_ib_i=e'$.

Suppose now that $p\geq\max\{e,e'\}$ and that $\alpha_p$ is the
image of $\alpha$ in $J_p(X)$. If we tensor everything with
$k[t]/(t^{p+1})$, we get the following factorization of the
pull-back map $h_p\colon
\alpha_p^*f^*\Omega_X\to\alpha_p^*\Omega_{X'}$
\begin{equation}\label{eq_cv10}
\alpha_p^*f^*\Omega_X\overset{g'_p}\to
T_p=T\otimes_{k\llbracket t\rrbracket}k[t]/(t^{p+1})\overset{g_p}\to
\alpha_p^*\Omega_{X'},
\end{equation}
with $g'_p$ surjective
 and ${\rm Ker}(g'_p)={\rm Ker}(h)\otimes
_{k\llbracket t\rrbracket}k[t]/(t^{p+1})\simeq \oplus_i k[t]/(t^{b_i})$.

 The following lemma will be needed in
the proof of Theorem~\ref{change_of_variable}.

\begin{lemma}\label{lemma_change_of_variable}
Let $f\colon X'\to X$ be as in the theorem. Suppose that $\gamma_m$,
$\gamma'_m\in J_m(X')$ are such that $\ord_{\gamma_m}({\rm
Jac}_f)=e$ and $\ord_{f_m(\gamma_m)}({\rm Jac}_X)=e'$, with
$m\geq\max\{2e,e+e'\}$. If $f_m(\gamma_m)=f_m(\gamma'_m)$, then
$\gamma_m$ and $\gamma'_m$ have the same image in $J_{m-e}(X')$.
\end{lemma}

\begin{proof}
For an arc $\delta$ we will denote by $\delta_m$ its image in the
space of $m$--jets.
 It is enough to show the following claim: if $q\geq \max\{2e,e+e'\}$, and if we have
 $$\alpha\in J_{\infty}(X'),\,\,\beta\in J_{\infty}(X),$$
with $\ord_{\alpha}({\rm Jac}_f)=e$, $\ord_{\beta}({\rm Jac}_X)=e'$
and $f_q(\alpha_q)=\beta_q$, then there is $\delta\in
J_{\infty}(X')$ having the same image as $\alpha$ in $J_{q-e}(X')$
and such that $f_{q+1}(\delta_{q+1})=\beta_{q+1}$.

Indeed, in the situation in the lemma, let us choose arbitrary
liftings $\gamma$ and $\gamma'$ of $\gamma_m$ and $\gamma'_m$,
respectively, to $J_{\infty}(X')$. We use the above claim to
construct recursively $\alpha^{(q)}\in J_{\infty}(X')$ for $q\geq m$
such that $\alpha^{(m)}=\gamma$ and $\alpha^{(q+1)}$, $\alpha^{(q)}$
have the same image in $J_{q-e}(X')$ and
$$f_q(\alpha^{(q)}_q)=\psi^X_q(f_{\infty}(\gamma'))$$
for every $q\geq m$ (note that since $m\geq\max\{2e,e+e'\}$ each
$\alpha^{(q)}$ vanishes along ${\rm Jac}_f$ and $f^{-1}({\rm
Jac}_X)$ with the same order as $\gamma$). The sequence given by the
image of each $\alpha^{(q)}$ in $J_q(X')$ defines a unique
$\alpha\in J_{\infty}(X')$ such that $\alpha$ and $\alpha^{(q)}$
have the same image in $J_{q-e}(X')$ for every $q\geq m$. We deduce
that $f_{\infty}(\alpha)=f_{\infty}(\gamma')$. Since $\alpha$ has
the same image as $\gamma$ in $J_{m-e}(X')$, and since $m-e\geq
\max\{e,e'\}$, it follows that
$$\alpha\not\in J_{\infty}(V({\rm Jac}_f))\cup
f_{\infty}^{-1}(J_{\infty}(X_{\rm sing})),$$ hence $\alpha=\gamma'$
by Lemma~\ref{injectivity2}. In particular, $\gamma$ and $\gamma'$
have the same image in $J_{m-e}(X')$.

We now prove  the claim made at the beginning of the proof. It
follows from Proposition~\ref{fiber6} i) that using $\alpha_{q+1}\in
(\pi^{X'}_{q+1,q-e})^{-1}(\alpha_{q-e})$ we get an isomorphism
$$(\pi^{X'}_{q+1,q-e})^{-1}(\alpha_{q-e})\simeq {\rm Hom}_{k[t]/(t^{q-e+1})}(\alpha_{q-e}^*\Omega_{X'},
(t^{q-e+1})/(t^{q+2})).$$ Similarly, using $f_{q+1}(\alpha_{q+1})$
we see that
$$(\pi^X_{q+1,q-e})^{-1}(\beta_{q-e}))\simeq
{\rm Hom}_{k[t]/(t^{q-e+1})}(\beta_{q-e}^*\Omega_X,
(t^{q-e+1})/(t^{q+2})).$$

Via this isomorphism $\beta_{q+1}$ corresponds to $w
\colon\beta_{q-e}^*\Omega_X\to (t^{q-e+1})/(t^{q+2})$. Note that
since $\beta_q=f_q(\alpha_q)$, the image of $w$ lies in
$(t^{q+1})/(t^{q+2})$. We now use  the factorization (\ref{eq_cv10})
with $p=q-e$. If we construct a morphism $u\colon
\alpha_{q-e}^*\Omega_{X'}\to (t^{q-e+1})/(t^{q+2})$ such that
$u\circ h_{q-e}=w$, then $u$ corresponds to an element
$\delta_{q+1}\in J_{q+1}(X')$ such that any lifting $\delta$ of
$\delta_{q+1}$ to $J_{\infty}(X')$ satisfies our requirement.

We first show that $w$ is zero on ${\rm Ker}(g'_{q-e})$. Note 
that by using $f_{2q+1}(\alpha)\in
(\pi^X_{2q+1,q})^{-1}(f_q(\alpha_q))$ we see that $\beta_{2q+1}$
corresponds to a morphism
$$w'\colon\beta_q^*\Omega_X\to (t^{q+1})/(t^{2q+2}),$$
such that $w$ is obtained by tensoring $w'$ with $k[t]/(t^{q-e+1})$
and composing with the natural map $(t^{q+1})/(t^{2q-e+2})\to
(t^{q-e+1})/(t^{q+2})$. Therefore in order to show that $w$ is zero
on ${\rm Ker}(g'_{q-e})$ it is enough to show that $w'$ maps ${\rm
Ker}(g'_{q-e})$ to $(t^{q+2})/(t^{2q+2})$. Since ${\rm
Ker}(g'_{q-e})$ is a direct sum of $k[t]/(t^{q+1})$--modules of the
form $k[t]/(t^b)$ with $b\leq e'$, it follows that $w'({\rm
Ker}(h'))$ is contained in $(t^{2q+2-e'})/(t^{2q+2})$. We have
$2q+2-e'\geq q+2$, hence $w$ is zero on ${\rm Ker}(g'_{q-e})$.

Therefore $w$ induces a morphism $\overline{w}\colon T_{q-e}\to
(t^{q+1})/(t^{q+2})$. We know that in suitable bases of $T_{q-e}$ and
$\beta_{q-e}^*\Omega_X$ the map $g_{q-e}$ is given by the diagonal
matrix with entries $t^{a_1},\ldots,t^{a_n}$, with all $a_i\leq e$.
It follows that we can find $u\colon\alpha_{q-e}^*\Omega_{X'}\to
(t^{q+1-e})/(t^{q+2})$ such that $u\circ h_{q-e}=w$, which completes the
proof of the lemma.
\end{proof}

\begin{proof}[Proof of Theorem~\ref{change_of_variable}]
The assertion in i) follows from
Lemma~\ref{lemma_change_of_variable}, and we now prove ii). We first show
that every fiber of the restriction of $f_m$ to
$\psi_m^{X'}(C_{e,e'})$ is isomorphic to $\AAA^e$, and we then
explain how to globalize the argument. Note first that since $X'$ is
nonsingular, every jet in $J_m(X')$ can be lifted to
$J_{\infty}(X')$, hence an element in $J_m(X')$ lies in
$\psi_m^{X'}(C_{e,e'})$ if and only if its projection to
$J_{m-e}(X')$ lies in $\psi_{m-e}^{X'}(C_{e,e'})$.

Let $\gamma'_m\in \psi_m^{X'}(C_{e,e'})$ and $\gamma'_{m-e}$ its
image in $J_{m-e}(X')$. We denote by $\gamma_m$ and $\gamma_{m-e}$
the images of $\gamma'_m$ and $\gamma'_{m-e}$ by $f_m$ and
$f_{m-e}$, respectively. It follows from
Lemma~\ref{lemma_change_of_variable} that $f_m^{-1}(\gamma_m)$ is
contained in the fiber of $\pi_{m,m-e}^{X'}$ over $\gamma'_{m-e}$.
Using the identifications of the fibers of $\pi_{m,m-e}^{X'}$ and
$\pi_{m,m-e}^X$ over $\gamma'_{m-e}$ and, respectively,
$\gamma_{m-e}$ given by Proposition~\ref{fiber6}, we get an
isomorphism of $f_m^{-1}(\gamma_m)$ with the kernel of
\begin{equation}\label{fiber_identification}
{\rm Hom}((\gamma'_{m-e})^*\Omega_{X'},(t^{m-e+1})/(t^{m+1}))\to
{\rm Hom}(\gamma_{m-e}^*\Omega_{X},(t^{m-e+1})/(t^{m+1})),
\end{equation}
where the Hom groups are over $k[t]/(t^{m-e+1})$. This gives an
isomorphism
$$f_m^{-1}(\gamma_m)\simeq {\rm
Hom}((\gamma'_{m-e})^*\Omega_{X'/X},(t^{m-e+1})/(t^{m+1})).$$ Since
$(\gamma'_{m-e})^*\Omega_{X'/X}\simeq
k[t]/(t^{a_1})\oplus\ldots\oplus k[t]/(t^{a_n})$, with $0\leq
a_i\leq\ldots\leq a_n\leq e$, with $\sum_ia_i=e$, we deduce
$$f_m^{-1}(\gamma_m)\simeq
\oplus_{i=1}^n(t^{m+1-a_i})/(t^{m+1})\simeq \AAA^e.$$

We now show that the above argument globalizes to give the full
assertion in ii). Note first that after restricting to an affine
open subset of $X'$, we may assume that we have a section of
$\pi_{m,m-e}^{X'}$. By Remark~\ref{remark_nonsingular_case}, it
follows that $J_m(X')$ becomes isomorphic to a geometric vector
bundle $E$ over $J_{m-e}(X')$ whose fiber over some $\gamma'_{m-e}$
is isomorphic to ${\rm
Hom}((\gamma'_{m-e})^*\Omega_{X'},(t^{m-e+1})/(t^{m+1}))$. Moreover,
after restricting to a suitable locally closed cover of
$\psi_{m-e}^{X'}(C_{e,e'})$, we may assume that, in the above
notation, the integers $a_1,\ldots,a_n$ do not depend  on
$\gamma_{m-e}'$. It follows that we get a geometric subbundle $F$ of
$E$ over this subset of $J_{m-e}(X')$ whose fiber over
$\gamma_{m-e}'$ is ${\rm
Hom}((\gamma'_{m-e})^*\Omega_{X'/X},(t^{m-e+1})/(t^{m+1}))$. It
follows from the above discussion that we get a one-to-one map from
the quotient bundle $E/F$ to $J_m(X)$. This completes the proof of
the theorem.
\end{proof}

\begin{corollary}\label{cor1_change_of_variable}
Suppose that $k$ is uncountable. With the notation in
Theorem~\ref{change_of_variable}, if $A\subseteq C_{e,e'}$ is a
cylinder in $J_{\infty}(X')$, then $f_{\infty}(A)$ is a cylinder in
$J_{\infty}(X)$.
\end{corollary}

\begin{proof}
Suppose that $A=\psi_p^{-1}(S)$, and let $m\geq\max\{2e,e+e',e+p\}$.
It is enough to show that
$$f_{\infty}(A)=(\psi_m^X)^{-1}(f_m(\psi_m^{X'}(A))).$$
The inclusion "$\subseteq$" is trivial, hence it is enough to show
the reverse inclusion. Consider $\delta\in
(\psi_m^X)^{-1}(f_m(\psi_m^{X'}(A)))$. In particular $\delta\not\in
J_{\infty}(X_{\rm sing})$, and by Corollary~\ref{cor2_cylinder1}
there is $\gamma\in J_{\infty}(X')$ such that
$\delta=f_{\infty}(\gamma)$. Since $f_m(\psi_m^{X'}(\gamma))\in
f_m(\psi_m^{X'}(A))$, it follows from
Lemma~\ref{lemma_change_of_variable} that the image of $\gamma$ in
$J_p(X')$ lies in $S$, hence $\gamma\in A$.
\end{proof}

\begin{corollary}\label{cor2_change_of_variable}
Suppose that $k$ is uncountable and of characteristic zero. With the
notation in the theorem, if $B\subseteq J_{\infty}(X)$ is a
cylinder, then
$$\codim(B)=\min\{\codim(f_{\infty}^{-1}(B)\cap C_{e,e'})+e\vert
e,e'\in\NN\}.$$ Moreover, we have
$$|B|=\sum_{e,e'}|f_{\infty}^{-1}(B)\cap C_{e'e'}|,$$
where the sum is over those $e$, $e'\in\NN$ such that
$\codim(f_{\infty}^{-1}(B)\cap C_{e,e'})+e=\codim(B)$.
\end{corollary}

\begin{proof}
It follows from the previous corollary that each $B\cap
f_{\infty}(C_{e,e'})$ is a cylinder and Lemma~\ref{injectivity2}
implies that these cylinders are disjoint. Moreover, the complement
in $B$ of their union is thin, so $\codim(B)=\min_{e,e'}
\codim(B\cap f_{\infty}(C_{e,e'}))$ by
Proposition~\ref{countable_union} and $|B|=\sum_{e,e'}|B\cap
f_{\infty}(C_{e,e'})|$, the sum being over those $e$ and $e'$ such
that $\codim(f_{\infty}(C_{e,e'})\cap B)=\codim(B)$. The fact that
$$\codim(f_{\infty}(C_{e,e'})\cap B)=\codim(C_{e,e'}\cap
f_{\infty}^{-1}(B))+e,\,\,|f_{\infty}(C_{e,e'})\cap B|=
|C_{e,e'}\cap f_{\infty}^{-1}(B)|$$ is a direct consequence of
Theorem~\ref{change_of_variable}.
\end{proof}

\begin{remark}
Note that we needed to assume  ${\rm char}(k)=0$ simply because
we used existence of resolution of singularities in proving the
basic properties of codimension of cylinders.
\end{remark}

\section{Minimal log discrepancies via arcs}

{}From now on we assume that the characteristic of the ground field
is zero, as we will make systematic use of existence of resolution
of singularities. We start by recalling some basic definitions in
the theory of singularities of pairs.

We work with pairs $(X,Y)$, where $X$ is a normal $\QQ$--Gorenstein
$n$--dimensional variety and $Y=\sum_{i=1}^sq_iY_i$  is a formal
combination with real numbers $q_i$ and proper closed subschemes
$Y_i$ of $X$. An important special case is when $Y$ is an
$\RR$-\emph{Cartier divisor}, i.e. when all $Y_i$ are defined by
locally principal ideals.
 We say that a pair $(X,Y)$ is effective if all $q_i$ are
nonnegative. Since $X$ is normal, we have a Weil divisor $K_X$ on
$X$, uniquely defined up to linear equivalence, such that $\OO(K_X)
\simeq i_*\Omega_{X_{\rm reg}}^n$, where $i\colon X_{\rm
reg}\hookrightarrow X$ is the inclusion of the smooth locus.
Moreover, since $X$ is $\QQ$--Gorenstein, we may and will fix a
positive integer $r$ such that $rK_X$ is a Cartier divisor.

Invariants of the singularities of $(X,Y)$ are defined using
\emph{divisors over X}: these are prime divisors $E\subset X'$,
where $f\colon X'\to X$ is a birational morphism and $X'$ is normal.
Every such divisor $E$ gives a discrete valuation ${\rm ord}_E$ of
the function field $K(X')=K(X)$, corresponding to the DVR
$\OO_{X',E}$. We identify two divisors over $X$ if they give the
same valuation of $K(X)$. In particular, we may always assume that
$X'$ and $E$ are both smooth. The \emph{center} of $E$ is the
closure of $f(E)$ in $X$ and it is denoted by $c_X(E)$.

Let $E$ be a divisor over $X$. If $Z$ is a closed subscheme of $X$,
then we define $\ord_E(Z)$ as follows: we may assume that $E$ is a
divisor on $X'$ and that the scheme-theoretic inverse image
$f^{-1}(Z)$ is a divisor. Then $\ord_E(Z)$ is the coefficient of $E$
in $f^{-1}(Z)$. If $(X,Y)$ is a pair as above, then we put
$\ord_E(Y):=\sum_iq_i\ord_{E}(Y_i)$. We also define
$\ord_E(K_{-/X})$ as the coefficient of $E$ in $K_{X'/X}$. Recall
that $K_{X'/X}$ is the unique $\QQ$--divisor supported on the
exceptional locus of $f$ such that $rK_{X'/X}$ is linearly
equivalent with $rK_{X'}-f^*(rK_X)$. Note that both $\ord_E(Y)$ and
$\ord_E(K_{-/X})$  do not depend on the particular $X'$ we have
chosen.

Suppose now that $(X,Y)$ is a pair and that $E$ is a divisor over
$X$. The \emph{log discrepancy} of $(X,Y)$ with respect to $E$ is
$$a(E;X,Y):=\ord_E(K_{-/X})-\ord_E(Y)+1.$$
If $W$ is a closed subset of $X$, and $\dim(X)\geq 2$, then the \emph{minimal log
discrepancy} of $(X,Y)$ along $W$ is defined by
$${\rm mld}(W;X,Y):=\inf\{a(E;X,Y)\mid E\,{\rm divisor}\,{\rm over} X,\,c_X(E)\subseteq W\}.$$
When $\dim(X)=1$ we use the same definition of minimal log discrepancy, unless the infimum
is negative, in which case we make the convention that $\mld(W; X,Y)=-\infty$
(see below for motivation).  
There are also other versions of minimal log discrepancies (see
\cite{ambro}), but the study of all these variants can be reduced to
the study of the above one. In what follows we give a quick
introduction to minimal log discrepancies, and refer for proofs and
details to \emph{loc. cit.}

\begin{remark}\label{remark_integral_closure}
If $\widetilde{Y}:=\sum_iq_i\widetilde{Y_i}$, where each
$\widetilde{Y_i}$ is defined by the integral closure of the ideal
defining $Y_i$, then $\ord_E(Y)=\ord_E(\widetilde{Y})$ for every
divisor $E$ over $X$. For basic facts about integral closure of
ideals, see for example \cite{lazarsfeld}, \S 9.6.A. We deduce that
we have ${\rm mld}(W;X,Y)=\mld(W;X,\widetilde{Y})$.
\end{remark}

It is an easy computation to show that if $E$ and $F$ are divisors
with simple normal crossings on $X'$ above $X$, and if $F_1$ is the
exceptional divisor of the blowing-up of $X'$ along $E\cap F$ (we
assume that this is nonempty and connected), then
$$a(F_1;X,Y)=a(E;X,Y)+a(F;X,Y).$$
We may repeat this procedure, blowing-up along the intersection of
$F_1$ with the proper transform of $E$. In this way we get divisors
$F_m$ over $X$ for every $m\geq 1$ with
$$a(F_m;X,Y)=m\cdot a(E;X,Y)+a(F;X,Y).$$
In particular, this computation shows that if $\dim(X)\geq 2$ and
$\mld(W;X,Y)<0$, then $\mld(W;X,Y)=-\infty$ (which explains our convention 
in the one-dimensional case).

A pair $(X,Y)$ is \emph{log canonical} (\emph{Kawamata log
terminal}, or klt for short) if and only if $\mld(X;X,Y)\geq 0$
(respectively, $\mld(X;X,Y)>0$). Note that for a closed subset $W$,
if $\mld(W;X,Y)\geq 0$ then for every divisor $E$ over $X$ such that
$c_X(E)\cap W\neq\emptyset$ we have $a(E;X,Y)\geq 0$. Indeed, if
this is not the case, then we can find a divisor $F$ on some $X'$
with $c_X(F)\subseteq W$ and such that $E$ and $F$ have simple
normal crossings and nonempty intersection. As above, we produce a
sequence of divisors $F_m$ with $c_X(F_m)\subseteq W$ and
$\lim_{m\to\infty}a(F_m;X,Y)=-\infty$.

This assertion can be used to show that $\mld(W;X,Y)\geq 0$ if and
only if there is an open subset $U$ of $X$ containing $W$ such that
$(U,Y\vert_U)$ is log canonical. In fact, we have the following more
precise proposition that allows computing minimal log discrepancies
via log resolutions.

\begin{proposition}\label{compute_via_resolutions}
Let $(X,Y)$ be a pair as above and $W\subseteq X$ a closed subset.
Suppose that $f\colon X'\to X$ is a proper birational morphism with
$X'$ nonsingular, and such that the union of $\cup_if^{-1}(Y_i)$, of
the exceptional locus of $f$ and of $f^{-1}(W)$ (in case $W\neq X$)
is a divisor with simple normal crossings. Write
$$f^{-1}(Y):=\sum_iq_if^{-1}(Y_i)=\sum_{j=1}^{d}\alpha_jE_j,\,\,K_{X'/X}
=\sum_{j=1}^d\kappa_j E_j.$$ For a nonnegative real number $\tau$,
we have $\mld(W;X,Y)\geq \tau$ if and only if the following
conditions hold:
\begin{enumerate}
\item For every $j$ such that $f(E_j)\cap W\neq\emptyset$ we have
$\kappa_j+1-\alpha_j\geq 0$.
\item For every $j$ such that $f(E_j)\subseteq W$ we have
$\kappa_j+1-\alpha_j\geq\tau$.
\end{enumerate}
\end{proposition}

\bigskip

We now turn to the description of minimal log discrepancies in terms
of codimensions of contact loci from \cite{EMY}. We assume that $k$
is uncountable. If $(X,Y)$ is a pair with $Y=\sum_{i=1}^sq_i Y_i$
and if $w=(w_i)\in\NN^s$, then we put ${\rm Cont}^{\geq
w}(Y):=\cap_i{\rm Cont}^{\geq w_i}(Y_i)$, which is clearly a
cylinder. We similarly define ${\rm Cont}^w(Y)$, ${\rm Cont}^w(Y)_m$
and ${\rm Cont}^{\geq w}(Y)_m$.

Recall that $rK_X$ is a Cartier divisor. We have a canonical map
$$\eta_r\colon (\Omega_X^n)^{\otimes r}\to
\OO(rK_X)=i_*((\Omega_{X_{\rm reg}}^n)^{\otimes r}).$$ 
We can write ${\rm Im}(\eta_r)=I_{Z_r}\otimes\OO(rK_X)$,
and the subscheme $Z_r$ defined by $I_{Z_r}$ is the
$r^{\rm th}$ \emph{Nash subscheme} of $X$.
It is clear that $I_{Z_{rs}}=I_{Z_r}^s$ for every $s\geq 1$.

Suppose that $W$ is a
proper closed subset of $X$, and let $f\colon X'\to X$ be a
resolution of singularities as in
Proposition~\ref{compute_via_resolutions} such that, in addition,
$f^{-1}(V({\rm Jac}_X))$ and $f^{-1}(Z_r)$ are divisors, having
simple normal crossings with the exceptional locus of $f$, with
$f^{-1}(Y)$ and with $f^{-1}(W)$.

\begin{lemma}\label{lem_EMY}{\rm (}\cite{EMY}{\rm )}
Let $(X,Y)$ be a pair and $f\colon X'\to X$ a resolution as above.
Write
$$f^{-1}(Y_i)=\sum_{j=1}^d\alpha_{i,j}E_j,\,\,K_{X'/X}=\sum_{j=1}^d
\kappa_jE_j,\,\,f^{-1}(Z_r)=\sum_{j=1}^d z_jE_j.$$ For every
$w=(w_i)\in\NN^s$ and $\ell\in\NN$ we have
$$\codim({\rm Cont}^{w}(Y)\cap {\rm Cont}^{\ell}(Z_r)\cap {\rm
Cont}^{\geq
1}(W))=\frac{\ell}{r}+\min_{\nu}\sum_j(\kappa_j+1)\nu_j,$$ where the
minimum is over those $\nu=(\nu_j)\in\NN^d$ with
$\sum_j\alpha_{i,j}\nu_j=w_i$ for all $i$ and $\sum_jz_j\nu_j=\ell$,
and such that $\cap_{\nu_j\geq 1}E_j\neq\emptyset$ and $\nu_j\geq 1$
for at least one $j$ with $f(E_j)\subseteq W$.
\end{lemma}

\begin{proof}
For every $\nu=(\nu_j)\in\NN^d$ we put ${\rm
Cont}^{\nu}(E)=\cap_j{\rm Cont}^{\nu_j}(E_j)$. Since $\sum_jE_j$ has
simple normal crossings, we see that ${\rm Cont}^{\nu}(E)$ is
nonempty if and only if $\cap_{\nu_j\geq 1}E_j\neq\emptyset$, and in
this case all irreducible components of ${\rm Cont}^{\nu}(E)$ have
codimension $\sum_j\nu_j$. Indeed, by Lemma~\ref{lem2} it is enough
to check this when $X=\AAA^n$ and the $E_j$ are coordinate
hyperplanes, in which case the assertion is clear.

Suppose that $\gamma\in {\rm Cont}^{\nu}(E)$, hence
$\ord_{f_{\infty}(\gamma)}(Y_i)=\sum_j\alpha_{i,j}\nu_j$ and
$\ord_{f_{\infty}(\gamma)}(Z_r)=\sum_jz_j\nu_j$. It is clear that
$f_{\infty}(\gamma)\in {\rm Cont}^{\geq 1}(W)$ if and only if there
is $j$ such that $\nu_j\geq 1$ and $E_j\subseteq f^{-1}(W)$.

By the definition of $Z_r$ we have ${\rm
Jac}_f^r=f^{-1}(I_{Z_r})\cdot\OO(-rK_{X'/X})$, hence
$${\rm ord}_{\gamma}({\rm
Jac}_f)=\sum_j\left(\frac{z_j}{r}+\kappa_j\right)\nu_j.$$ Moreover,
by our assumption, the order of vanishing of arcs in ${\rm
Cont}^{\nu}(E)$ along $f^{-1}V({\rm Jac}_X)$ is finite and constant.
It follows from Corollary~\ref{cor1_change_of_variable} and
Theorem~\ref{change_of_variable} that $f_{\infty}({\rm
Cont}^{\nu}(E))$ is a cylinder with
$$\codim\,f_{\infty}({\rm
Cont}^{\nu}(E))=\sum_j\frac{z_j}{r}\nu_j+\sum_j(\kappa_j+1)\nu_j.$$

By Lemma~\ref{injectivity2} the cylinders $f_{\infty}({\rm
Cont}^{\nu}(E))$ for various $\nu$ are mutually disjoint. If we take
the union over those $\nu$ such that $\sum_j\alpha_{i,j}\nu_j=w_i$
for all $i$ and $\sum_jz_j\nu_j=\ell$, with $\nu_j\geq 1$ for some
$E_j\subseteq f^{-1}(W)$, this union is contained in ${\rm
Cont}^w(Y)\cap {\rm Cont}^{\ell}(Z_r)\cap {\rm Cont}^{\geq 1}(W)$.
Moreover, its complement is contained in $\cup_jJ_{\infty}(f(E_j))$,
hence it is thin. The formula in the lemma now follows from
Proposition~\ref{countable_union}.
\end{proof}

\begin{theorem}\label{thm_EMY}{\rm (}\cite{EMY}{\rm )}
If $(X,Y)$ is a pair as above, and $W\subset X$ is a proper closed
subset, then
$$\mld(W;X,Y)=\inf_{w,\ell}\left\{\codim\left({\rm Cont}^w(Y)\cap {\rm
Cont}^{\ell}(Z_r)\cap {\rm Cont}^{\geq
1}(W)\right)-\frac{\ell}{r}-\sum_{i=1}^sq_iw_i\right\},$$ where the
minimum is over the $w=(w_i)\in\NN^s$ and $\ell\in\NN$. Moreover, if
this minimal log discrepancy is finite, then the infimum on the
right-hand side is a minimum.
\end{theorem}

\noindent If $X$ is nonsingular, then $Z_r=\emptyset$ and the description of
minimal log discrepancies in the theorem takes a particularly simple
form.

\begin{proof}[Proof of Theorem~\ref{thm_EMY}]
Let $f$ be a resolution as in Lemma~\ref{lem_EMY}. We keep the
notation in that lemma and its proof. We also put
$f^{-1}(Y)=\sum_j\alpha_jE_j$. Note that
$\alpha_j=\sum_i\alpha_{i,j}q_i$. After restricting to an open
neighborhood of $W$ we may assume that all $f(E_j)$ intersect $W$.

We first show  that $\mld(W;X,Y)$ is bounded above by the infimum in
the theorem. Of course, we may assume that $\mld(W;X,Y)$ is finite.
Therefore $\kappa_j+1-\alpha_j\geq {\rm mld}(W;X,Y)$ if
$f(E_j)\subseteq W$ and 
$\kappa_j+1-\alpha_j\geq 0$ for every $j$.

Let $\nu=(\nu_j)\in\NN^s$ be such that $\cap_{\nu_j\geq 1}E_j\neq\emptyset$,
and $\nu_j\geq 1$ for some $j$ with $f(E_j)\subseteq W$. In this
case we have
$$\sum_{j=1}^s(k_j+1)\nu_j\geq
\sum_j\alpha_j\nu_j+{\rm mld}(W;X,Y)\cdot\sum_{f(E_j)\subseteq
W}\nu_j\geq\sum_j\alpha_j\nu_j+{\rm mld}(W;X,Y).$$ 
If $\sum_j\alpha_{i,j}\nu_j=w_i$ for every $i$, and $\sum_jz_j\nu_j=\ell$, then
$\sum_j\alpha_j\nu_j=\sum_iq_iw_i$,  and
the formula in
 Lemma~\ref{lem_EMY} gives
$${\rm mld}(W;X,Y)\leq
\codim\left({\rm Cont}^w(Y)\cap {\rm Cont}^{\ell}(Z_r)\cap {\rm
Cont}^{\geq 1}(W)\right)-\sum_iq_iw_i-\frac{\ell}{r}.$$

Suppose now that we fix $E_j$ such that $f(E_j)\subseteq W$. If
$w_i=\alpha_{i,j}$ for every $i$ and if $\ell=z_j$, then it follows
from Lemma~\ref{lem_EMY} that
$$\codim({\rm Cont}^w(Y)\cap {\rm Cont}^{\ell}(Z_r)\cap{\rm
Cont}^{\geq 1}(W))\leq k_j+1+\frac{\ell}{r} $$
$$=\sum_iq_iw_i+\frac{\ell}{r}+a(E_j;X,Y).$$
Such an inequality holds for every divisor over $X$ whose center is
contained in $W$, and we deduce that if $\dim(X)\geq 2$, then
 the infimum in the theorem is $\leq {\rm mld}(W;X,Y)$
(note that the infimum does not depend on the particular resolution,
and every divisor with center in $W$ appears on some resolution).
Moreover, we see that if $a(E_j;X,Y)={\rm mld}(W;X,Y)$, then the
infimum is obtained for the above intersection of contact loci.

In order to complete the proof of the theorem, it is enough to show that
if $X$ is a curve, and if $a(W; X,Y)<0$, then the infimum in the theorem is $-\infty$. 
Note that in this case $W$ is a (smooth) point on $X$, and we may assume that
$Y_i=n_iW$ for some $n_i\in\ZZ$. Therefore our condition says that $\sum_iq_in_i>1$.
Since $\codim({\rm Cont}^m(W))=m$,  we see by taking $w_i=mn_i$ for all $i$ that
$$\codim\left({\rm Cont}^w(Y)\right)-\sum_iq_i w_i =m\left(1-\sum_iq_in_i\right)\to-\infty,$$
when we let $m$ go to infinity.
\end{proof}

\begin{remark}\label{remark_thm_EMY}
If $(X,Y)$ is an effective pair and $W\subset X$ is a proper closed
subset, then ${\rm mld}(W;X,Y)$ is equal to
 $$\inf_{w,\ell}\left\{\codim\left({\rm Cont}^{\geq w}(Y)\cap {\rm Cont}^{\ell}(Z_r)\cap
{\rm Cont}^{\geq
1}(W)\right)-\frac{\ell}{r}-\sum_{i=1}^sq_iw_i\right\},$$ where the
infimum is over all $w\in\NN^s$ and $\ell\in\NN$. Indeed, note that
we have
$$\bigsqcup_{w'}\left({\rm Cont}^{w'}(Y)\cap {\rm Cont}^{\ell}(Z_r)\cap {\rm
Cont}^{\geq 1}(W)\right)\subseteq {\rm Cont}^{\geq w}(Y)\cap {\rm
Cont}^{\ell}(Z_r)\cap {\rm Cont}^{\geq 1}(W),$$ where the disjoint
union is over $w'\in \NN^s$ such that $w'_i\geq w_i$ for every $i$.
Since the complement of this union is contained in
$\cup_iJ_{\infty}(Y_i)$, hence it is thin, our assertion follows
from Theorem~\ref{thm_EMY} via Proposition~\ref{countable_union}
(we also use the fact that since $(X,Y)$ is effective, if $w'_i\geq w_i$ for all $i$, then
$\sum_iq_iw'_i\geq\sum_iq_iw_i$).
\end{remark}

\begin{remark}
We have assumed in Theorem~\ref{thm_EMY} that $W$ is a proper closed
subset. In general, it is easy to reduce computing minimal log
discrepancies to this case, using the fact that if $X$ is
nonsingular and if $Y$ is empty, then $\mld(X;X,Y)=1$. Indeed, this
implies that if $(X,Y)$ is an arbitrary pair and if we take
$W=X_{\rm sing}\cup\bigcup_i Y_i$, then
$$\mld(X;X,Y)=\min\{\mld(W;X,Y),1\}.$$
For example, one can use this (or alternatively, one could just
follow the proof of Theorem~\ref{thm_EMY}) to show that the pair
$(X,Y)$ is log canonical if and only if for every $w\in\NN^s$ and
every $\ell\in\NN$, we have
$$\codim\left({\rm Cont}^w(Y)\cap {\rm Cont}^{\ell}(Z_r)\right)\geq
\frac{\ell}{r}+\sum_iq_iw_i.$$
\end{remark}

\begin{remark}{\rm (}\cite{EMY}{\rm )}
The usual set-up in Mori Theory is to work with a normal variety $X$
and a $\QQ$--divisor $D$ such that $K_X+D$ is $\QQ$--Cartier (see for
example \cite{kollar}). The results in this section have analogues
in that context. Suppose for simplicity that $D$ is effective, 
giving an embedding ${\mathcal O}_X\hookrightarrow {\mathcal O}_X(rD)$,
and that
$r(K_X+D)$ is Cartier. The image of the composition
$$(\Omega_X^n)^{\otimes r}\to(\Omega_X^n)^{\otimes r}\otimes {\mathcal O}_X(rD)\to
{\mathcal O}_X(r(K_X+D))$$
can be written an $I_T\otimes {\mathcal O}_X(r(K_X+D))$, for a closed subscheme $T$ of $X$. Arguing as above, one can then show that if $W$ is a proper closed subset of $X$, then
$${\rm mld}(W; X,D)=\inf_{e\in\NN}\left\{\codim\left({\rm Cont}^e(T)\cap {\rm Cont}^{\geq 1}(W)\right)-
\frac{e}{r}\right\}.$$
\end{remark}

\begin{example}
Suppose that $X$ is nonsingular and $Y$, $Y'$ are effective
combinations of closed subschemes of $X$. If $P$ is a point on $X$,
then
\begin{equation}\label{eq_example01}
{\rm mld}(P; X,Y+Y')\leq {\rm mld}(P; X,Y)+{\rm
mld}(P;X,Y')-\dim(X).
\end{equation}

\noindent Indeed, let us write $Y=\sum_iq_iY_i$ and $Y'=\sum_iq'_iY_i$, where
the $q_i$ and the $q'_i$ are nonnegative real numbers. If one of the minimal log discrepancies
on the right-hand side of (\ref{eq_example01}) is $-\infty$, then $\mld(P;X,Y+Y')=
-\infty$, as well. Otherwise, we can find
$w$ and $w'\in\NN^s$ and irreducible components
 $C$ of ${\rm Cont}^{\geq w}(Y)$ and $C'$ of ${\rm Cont}^{\geq
w'}(Y')$ such that $\codim(C)=\sum_iq_iw_i+{\rm mld}(P;X,Y)$ and
$\codim(C')=\sum_iq_iw'_i+{\rm mld}(P;X,Y')$. Note that $C\cap C'$
is nonempty, since it contains the constant arc over $P$. If $m\gg
0$, then $\psi_m(C\cap C')=\psi_m(C)\cap \psi_{m}(C')$, and using
the fact that the fiber $\pi_m^{-1}(P)$ of $J_m(X)$ over $P$ is
nonsingular, we deduce
$$\codim(\psi_m(C)\cap\psi_m(C'),\pi_m^{-1}(P))\leq
\codim(\psi_m(C),\pi_m^{-1}(P))+\codim(\psi_m(C'),\pi_m^{-1}(P)).$$
Since $C\cap C'\subseteq {\rm Cont}^{\geq w+w'}(Y+Y')$, we deduce
our assertion from Remark~\ref{remark_thm_EMY}.
\end{example}

\smallskip

Our next goal is to give a different interpretation of minimal log
discrepancies that is better suited for applications. The main
difference is that we replace cylinders in the space of arcs by
suitable locally closed subsets in the spaces of jets.

Recall that $Z_r$ is the $r^{\rm th}$ Nash subscheme of $X$. 
The  \emph{non-lci subscheme of $X$ of level $r$} is defined by the ideal
$J_r=(\overline{{\rm Jac}_X^r}\colon I_{Z_r})$, 
where we denote by $\overline{\fra}$ the integral closure of an ideal $\fra$.  It is shown in 
Corollary~\ref{cor3_appendix} in the Appendix that
$J_r\cdot I_{Z_r}$ and ${\rm Jac}_X^r$ have the same integral closure.
Note also that by Remark~\ref{char_lci}, the subscheme defined by $J_r$ is supported on the 
set of points $x\in X$ such that ${\mathcal O}_{X,x}$ is not locally complete intersection.
 It follows
from the basic properties of integral closure that given any ideal $\fra$, we have
$\ord_{\gamma}(\fra)=\ord_{\gamma}(\overline{\fra})$ for every arc
$\gamma\in J_{\infty}(X)$. In particular,
$\ord_{\gamma}(J_r)+\ord_{\gamma}(I_{Z_r})=r\cdot\ord_{\gamma}({\rm
Jac}_X)$.

\begin{theorem}\label{new_interpretation}
Let $(X,Y)$ be an effective pair and $r$ and $J_r$ as above. If $W$
is a proper closed subset of $X$, then

$$\mld(W;X,Y)=\inf\{(m+1)\dim(X)+\frac{e'}{r}-\sum_iq_iw_i$$
$$- \dim({\rm Cont}^{\geq w}(Y)_m\cap {\rm Cont}^e({\rm
Jac}_X)_m \cap {\rm Cont}^{e'}(J_r)_m\cap {\rm Cont}^{\geq
1}(W))_m\},$$ where the infimum is over those $w\in\NN^s$, and
$e,e',m\in\NN$ such that
 $m\geq\max\{2e,e+e',e+w_i\}$. Moreover, if this minimal log
 discrepancy is finite, then the above infimum is a minimum.
\end{theorem}

\noindent It will follow from the proof that the expression in the above
infimum does not depend on $m$, as long as
$m\geq\max\{2e,e+e',e+w_i\}$. Note also that $e$ comes up only in
the condition on $m$. The condition in the theorem simplifies when
$X$ is locally complete intersection, since $J_r=\OO_X$ by
Remark~\ref{char_lci} in the Appendix.

\begin{proof}[Proof of Theorem~\ref{new_interpretation}]
It follows from Theorem~\ref{thm_EMY} (see also
Remark~\ref{remark_thm_EMY}) that
$\mld(W;X,Y)=\inf_{w,\ell}\left\{\codim(C_{w,\ell})-\frac{\ell}{r}-\sum_iq_iw_i\right\}$,
where $w\in\NN^s$, $\ell\in\NN$, and
$$C_{w,\ell}={\rm Cont}^{\geq w}(Y)\cap {\rm Cont}^{\ell}(Z_r)\cap {\rm Cont}^{\geq 1}(W).$$
On the other hand, Proposition~\ref{countable_union} gives
$$\codim(C_{w,\ell})=\min_{e\in\NN}\codim(C_{w,\ell}\cap {\rm Cont}^{e}({\rm
Jac}_X)),$$ and for every $e$ we can write
$$C_{w,\ell}\cap {\rm Cont}^e({\rm Jac}_X)=
{\rm Cont}^{\geq w}(Y)\cap {\rm Cont}^e({\rm Jac}_X) \cap {\rm
Cont}^{e'}(J_r)\cap {\rm Cont}^{\geq 1}(W),$$ where $e'=re-\ell$.

Suppose now that $w$, $e$ and $\ell$ are fixed, $e'=re-\ell$, and
let $m\geq\max\{2e,e+e',e+w_i\}$. Consider  $$S:={\rm Cont}^{\geq
w}(Y)_m\cap {\rm Cont}^e({\rm Jac}_X)_m \cap {\rm
Cont}^{e'}(J_r)_m\cap {\rm Cont}^{\geq 1}(W)_m.$$
If we apply Proposition~\ref{fiber1} for the morphism
$\pi_{m,m-e}\colon J_m(X)\to J_{m-e}(X)$, we see that
$\dim(S)=\dim(\pi_{m,m-e}(S))+e(\dim(X)+1)$. Moreover,
$\pi_{m,m-e}(S)\subseteq {\rm Im}(\psi^X_{m-e})$ by
Proposition~\ref{fiber2}. It follows that
$$\codim(C_{w,\ell}\cap {\rm Cont}^e({\rm Jac}_X))=
(m-e+1)\dim(X)-\dim(\pi_{m,m-e}(S))$$ $$=(m+1)\dim(X) +e-\dim(S).$$
This gives the formula in the theorem.
\end{proof}

\begin{remark}\label{remark_new_interpretation}
If the pair $(X,Y)$ is not necessarily effective, then we can get an
analogue of Theorem~\ref{new_interpretation}, but involving contact
loci of specified order along each $Y_i$, as in
Theorem~\ref{thm_EMY}.
\end{remark}

\smallskip

In this section we have related the codimensions of various contact
loci with the numerical data of a log resolution. One can use, in
fact, Theorem~\ref{change_of_variable} to interpret also the "number
of irreducible components of minimal dimension" in the corresponding
contact loci. We illustrate this in the following examples. The
proofs are close in spirit to the proof of the other results in this section, so
we leave them for the reader.

\begin{example}\label{ex1_EMY}
Consider an effective pair $(X,Y)$ as above and $W\subset X$ a
proper closed subset. Suppose that $\tau:=\mld(W;X,Y)\geq 0$. We say
that a divisor $E$ over $X$ computes $\mld(W;X,Y)$ if
$c_X(E)\subseteq W$ and $a(E;X,Y)=\tau$. There is only one divisor
over $X$ computing $\mld(W;X,Y)$ if and only if for every
$w\in\NN^s$ and $m,e,e'\in\NN$ with $m\geq\max\{2e,e+e',e+w_i\}$,
there is at most one irreducible component of $${\rm Cont}^{\geq
w}(Y)_m\cap {\rm Cont}^e({\rm Jac}_X)_m \cap {\rm
Cont}^{e'}(J_r)_m\cap {\rm Cont}^{\geq 1}(W)_m$$ of dimension
$(m+1)\dim(X)+\frac{e'}{r}-\tau-\sum_iq_iw_i$. A similar equivalence
holds when $W=X$ and ${\rm mld}(X;X,Y)=0$.
\end{example}

\begin{example}\label{ex2_EMY}{\rm (}\cite{mustata}{\rm )}
Let $X$ be a nonsingular variety, and $Y\subset X$ a closed
subvariety of codimension $c$, which is reduced and irreducible.
Since
$$\dim\,J_m(Y)\geq\dim\,J_m(Y_{\rm reg})=(m+1)\dim(Y)$$
for every $m$, it follows from Theorem~\ref{thm_EMY} that $\mld(X;
X,cY)\leq 0$, with equality if and only if
$\dim\,J_m(Y)=(m+1)\dim(Y)$ for every $m$. In fact, note that if
$X'$ is the blowing-up of $X$ along $Y$, and if $E$ is the component
of the exceptional divisor that dominates $Y$, then $a(E;X,c Y)=0$.

Suppose now that $(X,c Y)$ is log canonical. The assertion in the
previous example implies that $E$ is the unique divisor over $X$
with $a(E;X,c Y)=0$ if and only if for every $m$, the unique
irreducible component of $J_m(Y)$ of dimension $(m+1)\dim(Y)$ is
$\overline{J_m(Y_{\rm reg})}$.

Assume now that $Y$ is locally complete intersection. Since $J_m(Y)$
can be locally defined in $J_m(X)$ by $c(m+1)$ equations, it follows
that every irreducible component of $J_m(Y)$ has dimension at least
$(m+1)\dim(Y)$. Hence $(X,c Y)$ is log canonical if and only if
$J_m(Y)$ has pure dimension for every $m$. In addition, we deduce
from the above discussion that $J_m(Y)$ is irreducible for every $m$
if and only if $(X,c Y)$ is log canonical and $E$ is the only
divisor over $X$ such that $a(E;X,c Y)=0$. It is shown in
\cite{mustata} that this is equivalent with $Y$ having rational
singularities.
\end{example}

\begin{example}\label{ex3_EMY}
Let $(X,Y)$ be an effective log canonical pair that is
\emph{strictly log canonical}, that is $\mld(X;X,Y)=0$. A
\emph{center of non-klt singularities} is a closed subset of $X$ of
the form $c_X(F)$, where $F$ is a divisor over $X$ such that
$a(F;X,Y)=0$. One can show that an irreducible closed subset
$T\subset X$ is such a center if and only if there are $w\in\NN^s$,
and $e, e'\in\NN$ not all zero, such that for
$m\geq\max\{2e,e+e',e+w_i\}$, some irreducible component  of
$${\rm Cont}^{\geq w}(Y)_m\cap{\rm Cont}^e({\rm Jac}_X)_m\cap
{\rm Cont}^{e'}(J_r)_m$$ has dimension
$(m+1)\dim(X)+\frac{e'}{r}-\sum_iq_iw_i$ and dominates $T$.
\end{example}

\section{Inversion of Adjunction}

We apply the description of minimal log discrepancies from the
previous section to prove the following version of Inversion of
Adjunction. This result has been proved also by Kawakita in
\cite{kawakita1}.

\begin{theorem}\label{inv_adj}
Let $A$ be a nonsingular variety and $X\subset A$ a closed normal
subvariety of codimension $c$. Suppose that $W\subset X$ is a proper
closed subset and $Y=\sum_{i=1}^sq_iY_i$ where all $q_i\in\RR_+$ and the
$Y_i\subset A$ are closed subschemes not containing $X$ in their
support. If $r$ is a positive integer such that $rK_X$ is Cartier
and if $J_r$ is the ideal defining the non-lci subscheme of level $r$ of $X$, then
$${\rm mld}\left(W;X,\frac{1}{r}V(J_r)+Y\vert_X\right)={\rm mld}(W; A,cX+Y),$$
where $Y\vert_X:=\sum_iq_i(Y_i\cap X)$.
\end{theorem}

When $X$ is locally complete intersection, then $J_r=\OO_X$, and we
recover the result from \cite{EM} saying that ${\rm
mld}(W;X,Y\vert_X)={\rm mld}(W;A,cX+Y)$. It is shown in \emph{loc.
cit.} that this is equivalent with the following version of
Inversion of Adjunction for locally complete intersection varieties.

\begin{corollary}
Let $X$ be a normal locally complete intersection variety and
$H\subset X$ a normal Cartier divisor. If $W\subset H$ is a proper
closed subset, and if $Y=\sum_{i=1}^sq_iY_i$, where all $q_i\in\RR_+$
and $Y_i$ are closed subsets of $X$ not containing $H$ in their
support, then
$${\rm mld}(W;H,Y\vert_H)={\rm mld}(W;X,Y+H).$$
\end{corollary}

For motivation and applications of the general case of the Inversion
of Adjunctions Conjecture, we refer to \cite{K+}. For results in the
klt and the log canonical cases, see \cite{kollar} and
\cite{kawakita1}.

We start with two lemmas. Recall that for every scheme $X$ we have a
morphism $\Phi_{\infty}\colon\AAA^1\times J_{\infty}(X)\to
J_{\infty}(X)$ such that if $\gamma$ is an arc lying over $x\in X$,
then $\Phi_{\infty}(0,\gamma)$ is the constant arc over $x$.

\begin{lemma}\label{lem1_inv_adj}
Let $X$ be a reduced, pure-dimensional scheme and $C\subseteq
J_{\infty}(X)$ a nonempty cylinder. If $\Phi_{\infty}(\AAA^1\times
C)\subseteq C$, then $C\not\subseteq J_{\infty}(X_{\rm sing})$.
\end{lemma}

\begin{proof}
Write $C=(\psi_m^X)^{-1}(S)$, for some $S\subseteq J_m(X)$. Let
$\gamma\in C$ be an arc lying over $x\in X$. By hypothesis, the
constant $m$--jet $\gamma_m^x$ over $x$ lies in $S$. We take a
resolution of singularities $f\colon X'\to X$. It is enough to show
that $f_{\infty}^{-1}(C)$ is not contained in
$f_{\infty}^{-1}(J_{\infty}(X_{\rm sing}))$.

Let $x'$ be a point in $f^{-1}(x)$. The constant jet $\gamma_m^{x'}$
lies in $f_m^{-1}(S)$, hence $C':=(\psi_m^{X'})^{-1}(\gamma_m^{x'})$
is contained in $f_{\infty}^{-1}(C)$. On the other hand, $X'$ is
nonsingular, hence $C'$ is not contained in
$f_{\infty}^{-1}(J_{\infty}(X_{\rm sing}))=J_{\infty}(f^{-1}(X_{\rm
sing}))$ by Lemma~\ref{cylinder1}.
\end{proof}

We will apply this lemma as follows. We will consider a reduced and
irreducible variety $X$ embedded in a nonsingular variety $A$. In
$J_{\infty}(A)$ we will take a finite intersection of closed
cylinders of the form ${\rm Cont}^{\geq m}(Z)$. Such an intersection
is preserved by $\Phi_{\infty}$, and therefore so is each
irreducible component $\widetilde{C}$. The lemma then implies that
$C:=\widetilde{C}\cap J_{\infty}(X)$ is not contained in
$J_{\infty}(X_{\rm sing})$.

\begin{lemma}\label{lem2_inv_adj}
Let $A$ be a nonsingular variety and $M=H_1\cap\ldots\cap H_c$ a
codimension $c$ complete intersection in $A$. If $C$ is an
irreducible locally closed cylinder in $J_{\infty}(A)$ such that
$$C\subseteq\bigcap_{i=1}^c{\rm Cont}^{\geq d_i}(H_i),$$
and if there is $\gamma\in C\cap J_{\infty}(M)$ with
$\ord_{\gamma}({\rm Jac}_M)=e$, then
$$\codim(C\cap
J_{\infty}(M),J_{\infty}(M))\leq\codim(C,J_{\infty}(A))+e-\sum_{i=1}^c
d_i.$$
\end{lemma}

\begin{proof}
We may assume that $e$ is the smallest order of vanishing along
${\rm Jac}_M$ of an arc in $C\cap J_{\infty}(M)$.
 Let $m\geq \max\{2e,e+d_i\}$ be such that
$C=(\psi_{m-e}^A)^{-1}(S)$ for some irreducible locally closed
subset $S$ in $J_{m-e}(A)$. Let $S'$ be the inverse image of $S$ in
$J_m(A)$ and $S''$ an irreducible component of $S'\cap J_m(M)$
containing some jet having order $e$ along ${\rm Jac}_M$. Every jet
in $S'$ has order $\geq d_i$ along $H_i$, hence $S'\cap J_m(M)$ is
cut out in $S'$ by $\sum_i(m-d_i+1)$ equations, and therefore
$$\dim(S'')\geq \dim(S')-(m+1)c+\sum_{i=1}^c d_i=\dim(S)+e\cdot\dim(A)-(m+1)c+\sum_{i=1}^c d_i.$$

Let $S''_0$ be the open subset of $S''$ consisting of jets having
order $\leq e$ along ${\rm Jac}_M$. It follows from
Proposition~\ref{fiber2} (see also Remark~\ref{remark_fiber2}) that
the image in $J_{m-e}(M)$ of any element in $S''_0$ can be lifted to
$J_{\infty}(M)\cap C$, hence by assumption its order of vanishing
along ${\rm Jac}_M$ is $e$. Moreover, Proposition~\ref{fiber1}
implies that the image of $S''_0$ in $J_{m-e}(M)$ has dimension
$\dim(S''_0)-e(\dim(M)+1)$. We conclude that
$$\codim(C\cap J_{\infty}(M), J_{\infty}(M))\leq
(m-e+1)\dim(M)-\dim(S''_0)+e(\dim(M)+1)\leq$$
$$(m-e+1)\dim(A)+e-\dim(S)-\sum_{i=1}^c d_i=
\codim(C,J_{\infty}(A))+e-\sum_{i=1}^c d_i.$$
\end{proof}

\begin{proof}[Proof of Theorem~\ref{inv_adj}]
The assertion is local, hence we may assume that $A$ is affine. We
first show that ${\rm mld}(W; X,\frac{1}{r}V(J_r)+Y\vert_X) \geq
{\rm mld}(W; A,cX+Y)$. Suppose that this is not the case, and let us
use Theorem~\ref{new_interpretation} for
$(X,\frac{1}{r}V(J_r)+Y\vert_X)$. We get $w\in\NN^s$ and $e$, $e'$,
$m\in\NN$ such that $m\geq \max\{2e,e+e',e+w_i\}$ and $S\subseteq
J_m(X)$ with
$$S\subseteq {\rm Cont}^{\geq w}(Y)_m\cap {\rm Cont}^e({\rm Jac}_X)_m
\cap {\rm Cont}^{e'}(J_r)\cap {\rm Cont}^{\geq 1}(W)$$ such that
$\dim(S)>(m+1)\dim(X)-{\rm mld}(W;A,cX+Y)-\sum_iq_iw_i$. We may
consider $S$ as a subset of $J_m(A)$ contained in ${\rm Cont}^{\geq
(m+1)}(X)$, and applying Theorem~\ref{new_interpretation} for the pair $(A,cX+Y)$
we see that
$$\dim(S)\leq (m+1)\dim(A)-c(m+1)-\sum_iq_iw_i-{\rm
mld}(W;A,cX+Y).$$ This gives a contradiction.

We now prove the reverse inequality
 $$\tau:={\rm mld}(W; X,\frac{1}{r}V(J_r)+Y\vert_X) \leq
{\rm mld}(W; A,cX+Y).$$ If this does not hold, then we apply
Theorem~\ref{thm_EMY} (see also Remark~\ref{remark_thm_EMY}) to find
$w\in \NN^s$ and $d\in\NN$ such that for some irreducible component
$C$ of ${\rm Cont}^{\geq w}(Y)\cap {\rm Cont}^{\geq d}(X)$ we have
$\codim(C)<cd+\sum_iq_iw_i+\tau$. It follows from
Lemma~\ref{lem1_inv_adj} that $C\cap J_{\infty}(X)\not\subseteq
J_{\infty}(X_{\rm sing})$. Let $e$ be the smallest order of
vanishing along ${\rm Jac}_X$ of an arc in $C\cap J_{\infty}(X)$.
Fix such an arc $\gamma_0$.

Consider the closed subscheme $M\subset A$ whose ideal $I_M$ is
generated by $c$ general linear combinations of the generators of
the ideal $I_X$ of $X$. Therefore $M$ is a complete intersection and
$\ord_{\gamma_0}({\rm Jac}_M)=e$. By Corollary~\ref{cor1_appendix}
in the Appendix, we have ${\rm Jac}_M\cdot\OO_X\subseteq\left( (I_M\colon
I_X)+I_X\right)/I_X$. It follows that $\gamma_0$ lies in the cylinder
$$C_0:=C\cap {\rm Cont}^{\leq e}({\rm Jac}_M)\cap{\rm Cont}^{\leq e}(I_M\colon I_X).$$
$C_0$ is a nonempty open subcylinder of $C$, hence
$\codim(C)=\codim(C_0)$. On the other hand, Lemma~\ref{lem2_inv_adj}
gives
$$\codim(C_0\cap J_{\infty}(M), J_{\infty}(M))\leq
\codim(C_0)+e-cd.$$

If $\gamma\in J_{\infty}(M)$, then $\ord_{\gamma}({\rm Jac}_X)\leq
\ord_{\gamma}({\rm Jac}_M)$. If $\gamma$ lies also in $C_0$, then
$\gamma$ can't lie in the space of arcs of any other irreducible
component of $M$ but $X$ (we use the fact that $\gamma$ has finite
order along $(I_M\colon I_X)$, and the support of the scheme defined
by $(I_M\colon I_X)$ is the union of the irreducible components of
$M$ different from $X$). Therefore $C_0\cap J_{\infty}(M)= C_0\cap
J_{\infty}(X)$, and for every arc $\gamma$ in this intersection we
have $\ord_{\gamma}({\rm Jac}_X)=\ord_{\gamma}({\rm Jac}_M)=e$. We
deduce that
$$\codim(C_0\cap J_{\infty}(X),J_{\infty}(X))=\codim(C_0\cap
J_{\infty}(M),J_{\infty}(M))<\sum_iq_iw_i+\tau+e.$$

Since $C_0\cap J_{\infty}(X)=\cup_{e'=0}^{re}\left(C_0\cap
J_{\infty}(X)\cap {\rm Cont}^{e'}(J_r)\right)$, it follows that
there is $e'$ such that $\codim(C_0\cap J_{\infty}(X)\cap {\rm
Cont}^{e'}(J_r))<\sum_iq_iw_i+\tau+e$. On the other hand, this
cylinder is contained in ${\rm Cont}^{re-e'}(Z_r)$. We deduce from
Theorem~\ref{thm_EMY} (see also Remark~\ref{remark_thm_EMY}) that
${\rm mld}(W; X, \frac{1}{r}V(J_r)+Y\vert_X)<\tau$, a contradiction.
This completes the proof of the theorem.
\end{proof}

\begin{remark}
It follows from the above proof that even if the coefficients of $Y$
are negative, we still have the inequality
$${\rm mld}\left(W;X,\frac{1}{r}V(J_r)+Y\vert_X\right)\geq{\rm mld}(W; A,cX+Y)$$
(it is enough to use the description of minimal log discrepancies
mentioned in Remark~\ref{remark_new_interpretation}).
\end{remark}

\section{Appendix}

\subsection{Dimension of constructible subsets}

We recall here a few basic facts about the dimension of
constructible subsets. Let $X$ be a scheme of finite type over $k$,
and $W\subseteq X$ a constructible subset, with the induced Zariski
topology from $X$. If $A$ is a closed subset of $W$, we have
$\overline{A}\cap W=A$. Since $X$ is a Noetherian topological space
of bounded dimension, it follows that so is $W$.

Note that we have $\dim(W)=\dim(\overline{W})$. Indeed, the
inequality $\dim(W)\leq\dim(\overline{W})$ follows as above, while
the reverse inequality is a consequence of the fact that $W$
contains a subset $U$ that is open and dense in $\overline{W}$. We
see that if $W=T_1\cup\ldots\cup T_r$, where all $T_i$ are locally
closed (or more generally, constructible) in $X$, then
$\dim(W)=\max_i\{\dim(T_i)\}$.

Since $W$ is Noetherian, we have a unique decomposition
$W=W_1\cup\ldots \cup W_s$ in irreducible components. If $\dim(W)=n$
and if we have a decomposition $W=T_1\sqcup\ldots\sqcup T_r$ into disjoint
constructible subsets of $X$, then every irreducible component $A$
of some $T_i$, with $\dim(A)=n$ gives an irreducible component of
$W$ of dimension $n$, namely $\overline{A}\cap W$. Moreover, every
$n$--dimensional irreducible component of $W$  comes from a
unique $T_i$ and a unique such irreducible component of $T_i$.

If $f\colon X'\to X$ is a morphism of schemes that induces a
bijection between the constructible subsets $V'\subseteq X'$ and
$V\subseteq X$, then $\dim(V)=\dim(V')$ and $T\to V\cap
\overline{f(T)}$ gives a bijection between the irreducible
components of maximal dimension of $V'$ and those of $V$. It follows
that if we have
 a morphism of schemes $g\colon X'\to X$
and constructible subsets $V'\subseteq X'$ and $V\subseteq X$ such
that we get a weakly piecewise trivial fibration $V'\to V$ with
fiber $F$, then $\dim(V)=\dim(V')-\dim(F)$. Moreover, if $F$ is
irreducible, then we have a bijection between the irreducible
components of maximal dimension of $V'$ and those of $V$.

\subsection{Differentials and the canonical sheaf}

We start by reviewing the definition and some basic properties of
the canonical sheaf. The standard reference for this is
\cite{hartshorne}. To every pure-dimensional scheme over $k$ one
associates a coherent sheaf $\omega_X$ with the following
properties:

\begin{enumerate}
\item[i)] If $X$ is nonsingular of dimension $n$, then
there is a canonical isomorphism $\omega_X\simeq\Omega_X^n$.
\item[ii)] The definition is local: if $U$ is an open subset of $X$,
then there is a canonical isomorphism
$\omega_U\simeq\omega_X\vert_U$.
\item[iii)] If $X\hookrightarrow M$ is a closed subscheme of codimension $c$,
where $M$ is a pure-dimensional Cohen-Macaulay scheme, then there is
a canonical isomorphism
$$\omega_X\simeq{\mathcal Ext}^c_{\OO_M}(\OO_X,\omega_M).$$
\item[iv)] If $f\colon X\to M$ is a finite surjective morphism
of equidimensional schemes, then
$$f_*\omega_X\simeq {\mathcal Hom}_{\OO_M}(f_*\OO_X,\omega_M).$$
\item[v)] If $X$ is normal of dimension $n\geq 2$, then ${\rm depth}(\omega_X)\geq 2$.
Therefore there is a canonical isomorphism
$$\omega_X\simeq i_*\Omega_{X_{\rm reg}}^n,$$
where $i\colon X_{\rm reg}\hookrightarrow X$ is the inclusion of the
nonsingular locus of $X$.
\item[vi)] If $X$ is Gorenstein, then $\omega_X$ is locally free of
rank one.
\end{enumerate}

Note that $\omega_X$ is uniquely determined by properties i), ii)
and iii) above. Indeed, by ii) it is enough to describe
$\omega_{U_i}$ for the elements of an affine open cover $U_i$ of
$X$. On the other hand, if we embed $U_i$ as a closed subscheme of
codimension $c$ of an affine space $\AAA^{N}$, then we have
$$\omega_{U_i}\simeq {\mathcal
Ext}^c_{\OO_{\AAA^N}}(\OO_{U_i},\Omega_{\AAA^N}^N).$$ Note also that if
$Z$ is an irreducible component of $X$ that is generically reduced,
then by i) and ii) we see that the stalk of $\omega_X$ at the
generic point of $Z$ is $\Omega_{K/k}^n$, where $n=\dim(X)$ and $K$
is the residue field at the generic point of $Z$.

\bigskip

Suppose now that we are in the following setting. $X$ is a reduced
scheme of pure dimension $n$,
 and we have a closed embedding
$X\hookrightarrow A$, where $A$ is nonsingular of dimension $N$ and
has global algebraic coordinates $x_1,\ldots,x_N\in\Gamma(\OO_A)$
(that is, $dx_1,\ldots,dx_N$ trivialize $\Omega_A$). We assume that
the ideal $I_X$ of $X$ in $A$ is generated by
 $f_1,\ldots,f_d\in\Gamma(\OO_A)$.
For example, if $X$ is affine we may consider a closed embedding in
an affine space.

Let $c=N-n$. As in \S 4, for $1\leq i\leq d$ we  take
$F_i:=\sum_{j=1}^d a_{i,j}f_j$, where the $a_{i,j}$ are general
elements in $k$. If $M$ is the closed subscheme defined by
$I_M=(F_1,\ldots,F_c)$, then we have the following properties.

\begin{enumerate}
\item[1)] All irreducible components of $M$ have dimension $n$,
hence $M$ is a complete intersection.
\item[2)] $X$ is a closed subscheme of $M$ and $X=M$ at the generic
point of every irreducible component of $X$.
\item[3)] Some minor $\Delta$ of the Jacobian matrix of
$F_1,\ldots,F_c$ with respect to the coordinates $x_1,\ldots,x_N$
(let's say $\Delta={\rm det}(\partial F_i/\partial x_j)_{i,j\leq
c}$) does not vanish at the generic point of any irreducible
component of $X$.
\end{enumerate}

Moreover, every $c$ of the $F_i$ will satisfy similar properties.
Let us fix now $F_1,\ldots,F_c$ as above, generating the ideal
$I_M$. We also consider the residue scheme $X'$ of $X$ in $M$
defined by the ideal $(I_M\colon I_X)$. Note that $X'$ is supported
on the union of the irreducible components of $M$ that are not
contained in $X$. The intersection of $X$ and $X'$ is cut out in $X$
by the ideal $((I_M\colon I_X)+I_X)/I_X$.

Let $K$ denote the fraction field of $X$, i.e. $K$ is the product of
the residue fields of the generic points of the irreducible
components of $X$. We have a localization map
$\Omega_X^n\to\Omega_{K/k}^n$ given by taking a section of
$\Omega_X^n$ to its images in the corresponding stalks. By our
assumption $\Delta$ is an invertible element in $K$, and
$\Omega_{K/k}^n$ is freely generated over $K$ by
$dx_{c+1}\wedge\ldots\wedge dx_N$.

\begin{proposition}\label{prop1_appendix}
With the above notation, there are canonical morphisms
$$\Omega_X^n\overset{\eta}\to\omega_X\overset{u}\to\omega_M\vert_X
\overset{w}\to\Omega_{K/k}^n$$ with the following properties:
\begin{enumerate}
\item[a)] If $X$ is normal, then $\eta$
is given by the canonical isomorphism $\omega_X\simeq
i_*\Omega_{X_{\rm reg}}^n$, where $i\colon X_{\rm
reg}\hookrightarrow X$ is the inclusion of the nonsingular locus of
$X$.
\item[b)] $w$ is injective and identifies $\omega_M\vert_X$
with $\OO_X\cdot \Delta^{-1}dx_{c+1}\wedge\ldots\wedge dx_N$.
\item[c)] $u$ is injective and the image of $w\circ u$ is
$((I_M\colon I_X)+I_X)/I_X\cdot\Delta^{-1}dx_{c+1}\wedge\ldots\wedge
dx_N$.
\item[d)] The composition $w\circ u\circ\eta$ is the localization
map.  Its image is ${\rm Jac}(F_1,\ldots,F_c)\cdot\Delta^{-1}
dx_{c+1}\wedge\ldots\wedge dx_N$, where ${\rm Jac}(F_1,\ldots,F_c)$
denotes the ideal generated in $\OO_X$ by the $r$--minors of the
Jacobian matrix of $F_1,\ldots,F_c$.
\end{enumerate}
\end{proposition}

\begin{corollary}\label{cor1_appendix}
With the above notation, we have the following inclusion
$${\rm Jac}(F_1,\ldots,F_c)\cdot{\mathcal O}_X\subseteq ((I_M\colon I_X)+I_X)/I_X.$$
\end{corollary}

\begin{corollary}\label{cor2_appendix}
Suppose that $X$ is a normal affine $n$--dimensional Gorenstein variety. If
$Z$ is the first Nash subscheme of $X$, that is,  $I_Z\otimes\omega_X$ is
the image of the canonical map $\eta\colon\Omega_X^n\to\omega_X$,
then there is an ideal $J$ such that
$${\rm Jac}_X=I_Z\cdot J.$$
\end{corollary}

\begin{proof}
We choose a closed embedding $X\hookrightarrow A=\AAA^N$, and
let $F_1,\ldots,F_d$ be as above. For every $L= (i_1,\ldots,i_c)$,
with $1\leq i_1<\cdots<i_c\leq d$, let $I_L$ denote the ideal
generated by $F_{i_1},\ldots,F_{i_c}$. It follows from
Proposition~\ref{prop1_appendix} that
$${\rm Jac}(F_{i_1},\ldots,F_{i_c})\cdot {\mathcal O}_X=I_Z\cdot ((I_L\colon
I_X)+I_X)/I_X.$$ If we take $J=\sum_L((I_L\colon I_X)+I_X)/I_X$,
this ideal satisfies the condition in the corollary.
\end{proof}

\begin{proof}[Proof of Proposition~\ref{prop1_appendix}]
Since $X$ is reduced, we may consider its normalization $f\colon
\widetilde{X}\to X$. On $\widetilde{X}$ we have a canonical morphism
$\widetilde{\eta}\colon\Omega^n_{\widetilde{X}}\to\omega_{\widetilde{X}}$.
On the other hand, since $f$ is finite and surjective we have an
isomorphism $f_*\omega_{\widetilde{X}}\simeq {\mathcal Hom}_{\OO_X}
(f_*\OO_{\widetilde{X}},\omega_X)$, and the inclusion
$\OO_X\hookrightarrow f_*\OO_{\widetilde{X}}$ induces a morphism
$f_*\omega_{\widetilde{X}}\to \omega_X$.

 The morphism
$\eta$ is the composition
$$\Omega_X^n\to f_*\Omega_{\widetilde{X}}^n\overset{f_*\widetilde{\eta}}\to
f_*\omega_{\widetilde{X}}\to\omega_X,$$ where the first arrow is
induced by pulling-back differential forms. The construction is
compatible with the restriction to an open subset. In particular,
the composition
$$\Omega_X^n\to\omega_X\to\Omega^n_{K/k}$$
of $\eta$ with the morphism going to the stalks at the generic
points of the irreducible components of $X$ is the localization
morphism corresponding to $\Omega_X^n$.

Note that $\omega_M\simeq{\mathcal
Ext}^c_{\OO_A}(\OO_M,\Omega^N_A)$. Since $F_1,\ldots,F_c$ form a
regular sequence, we can compute $\omega_M$ using the Koszul complex
associated to the $F_i$'s to get
$$\omega_M\simeq {\mathcal
Hom}_{\OO_M}\left(\bigwedge^c(I_M/I_M^2),\Omega_A^N\vert_M\right).$$ This is a
free $\OO_M$-module generated by the morphism $\phi$ that takes
$\overline{F_1}\wedge\ldots\wedge\overline{F_c}$ to
$dx_1\wedge\ldots\wedge dx_N\vert_M$.

Since $X$ is a closed subscheme of $M$ of the same dimension, and
since $M$ is Cohen-Macaulay, it follows that $\omega_X\simeq
{\mathcal Hom}_{\OO_M}(\OO_X,\omega_M)$. In particular, we have
$\omega_X\subseteq\omega_M$. Moreover,
$$\omega_X\otimes\omega_M^{-1}
\simeq {\mathcal Hom}_{\OO_M}(\OO_X,\OO_M)=(I_M\colon I_X)/I_M.$$

Let $u$ be the composition $\omega_X\hookrightarrow\omega_M\to
\omega_M\vert_X$. Since $M=X$ at the generic point of each
irreducible component of $X$, $u$ is generically an isomorphism. On
the other hand, $M$ is Cohen-Macaulay and $\omega_X$ is contained in
the free $\OO_M$-module $\omega_M$, hence $\omega_X$ has no embedded
associated primes. Therefore $u$ is injective.

Using again the fact that $u$ is an isomorphism at the generic
points of the irreducible components of $X$ we get a localization
morphism $w\colon\omega_M\vert_X\to\Omega_{K/k}^n$, and we see that
the composition $w\circ u\circ\eta$ is the localization map for
$\Omega_X^n$ at the generic points.

By construction, $w$ takes the image of $\phi$ in $\omega_M\vert_X$
 to $\Delta^{-1}dx_{c+1}\wedge\ldots\wedge dx_N$. It follows from
our previous discussion that the image of $\omega_X$ in
$\omega_M\vert_X$ is $((I_M\colon
I_X)+I_X)/I_X\cdot\omega_M\vert_X$, from which we get the image of
$w\circ u$. The last assertion in d) follows from the fact that if
$1\leq i_1<\cdots<i_n\leq N$ and if $D$ is the $r$--minor of the
Jacobian of $F_1,\ldots,F_c$ corresponding to the variables
different from $x_{i_1},\ldots,x_{i_n}$, then
$$(w\circ u\circ\eta)(dx_{i_1}\wedge\ldots\wedge
dx_{i_n})=\pm \frac{D}{\Delta}dx_{c+1}\wedge\ldots\wedge dx_N.$$
This completes the proof of the proposition.
\end{proof}

\smallskip

Suppose now that $X$ is an affine $\QQ$--Gorenstein normal variety.
Our goal is to generalize Corollary~\ref{cor2_appendix} to this
setting. Let $K_X$ be a Weil divisor on $X$ such that
$\OO(K_X)\simeq\omega_X$ and let
 us fix a positive integer $r$ such that $r K_X$ is Cartier. Note
 that we have a canonical morphism $p_{r}\colon\omega_X^{\otimes r}
 \to\OO(r K_X)$.

We use the notation in Proposition~\ref{prop1_appendix}. Let
$\eta_{r} \colon(\Omega_X^n)^{\otimes r}\to \OO(r K_X)$ be the
composition of $\eta^{\otimes r}$ with $p_{r}$. Equivalently, if $i$
denotes the inclusion of $X_{\rm reg}$ into $X$, then $\eta_{r}$ is
identified with the canonical map $(\Omega_X^n)^{\otimes r}\to
i_*((\Omega_{X_{\rm reg}}^n)^{\otimes r})$. The image of $\eta_r$ is 
by definition $I_{Z_r}\otimes\OO(r K_X)$, where $Z_{r}$ is the $r^{\rm th}$ Nash subscheme of 
$X$.

Since $\omega_M^{\otimes r}\vert_X$ is locally-free, the morphism
$u^{\otimes r}$ induces 
$$u_r=i_*(u^{\otimes r}\vert_{X_{\rm reg}})
\colon \OO(rK_X)\to \omega_M^{\otimes r}\vert_X.$$
This is injective, since this is the case if we restrict
to the nonsingular locus of $X$. If we put $w_r:=w^{\otimes r}$,
then it follows from Proposition~\ref{prop1_appendix} that
\begin{enumerate}
\item[i)] $w_r$ is injective and its image is
$\OO_X\cdot\Delta^{-r}(dx_{c+1}\wedge\ldots\wedge dx_N)^{\otimes
r}$.
\item[ii)] The composition $w_r\circ u_r\circ\eta_r$ is the
localization map. Moreover, its image is equal to ${\rm
Jac}(F_1,\ldots,F_c)^r\cdot\Delta^{-r}(dx_{c+1}\wedge\ldots\wedge
dx_N)^{\otimes r}$.
\end{enumerate}

We now generalize  Corollary~\ref{cor2_appendix} to the case when $X$
is $\QQ$--Gorenstein. Let $\overline{\fra}$ denote the integral
closure of an ideal $\fra$. We define the \emph{non-lci subscheme of level $r$} to be 
the subscheme of $X$ defined by the ideal 
$J_r=(\overline{{\rm Jac}_X^r}\colon I_{Z_r})$ (see Remark~\ref{char_lci} below
for a justification of the name).

\begin{corollary}\label{cor3_appendix}
Let $X$ be a normal $\QQ$--Gorenstein $n$--dimensional variety and let
$r$ be a positive integer such that $rK_X$ is Cartier. If $Z_r$ is
the $r^{\rm th}$ Nash subscheme of $X$, and if $J_r$ defines the non-lci subscheme of level $r$, then 
the ideals ${\rm
Jac}_X^r$ and $I_{Z_r}\cdot J_r$ have the same integral closure.
\end{corollary}

\begin{proof}
It is enough to prove the assertion when $X$ is affine, hence we may
assume that we have a closed embedding  $X\subset A$ of codimension
$c$, and general elements $F_1,\ldots,F_d$ that generate the ideal
of $X$ in $A$, as above. It is enough to show that 
for every $L= (i_1,\ldots,i_c)$ with $1\leq i_1<\cdots<i_c\leq d$
we can find an ideal $\frb_L$ such that 
\begin{equation}\label{eq_cor33}
I_{Z_r}\cdot \frb_L=
{\rm Jac}(F_{i_1},\ldots,F_{i_c})^r.
\end{equation}
Indeed, in this case if we put $\frb:=\sum_L\frb_L$, then 
${\rm Jac}_X^r$  and $I_{Z_r}\cdot\frb$ have the same integral closure.
In particular, we have $\frb\subseteq J_r$, and we see that the inclusions
$$I_{Z_r}\cdot \fra_r\subseteq I_{Z_r}\cdot J_r\subseteq
\overline{{\rm Jac}_X^r}$$ 
become equalities after passing to integral closure. Note also that $\frb$ and $J_r$ have the
same integral closure.

In order to find $\frb_L$, we may assume without any loss of generality that
$L=(1,\ldots,c)$. With the above notation, consider the factorization
of the localization map $(\Omega_X^n)^{\otimes r}\to(\Omega_{K/k}^n)^{\otimes r}$
as $w_r\circ u_r\circ\eta_r$.
If $\frb_L$ is the ideal of $\OO_X$ such that the image of $u_r$ is $\frb_L\otimes\omega_M^{\otimes r}\vert_X$,
then (\ref{eq_cor33}) follows from the discussion preceding the statement of the corollary.
\end{proof}

\begin{remark}
Since
$I_{Z_{rs}}=I_{Z_r}^s$ for every $s\geq 1$, it follows that $(J_r)^s\subseteq J_{rs}$,
and we deduce from the corollary that these two ideals have the same
integral closure.
\end{remark}

\begin{remark}\label{char_lci}
Under the assumptions in Corollary~\ref{cor3_appendix}, the support
of the non-lci subscheme of level $r$
is the set of points $x\in X$ such
that $\OO_{X,x}$ is not locally complete intersection. Indeed, if
$\OO_{X,x}$ is locally complete intersection, then after replacing
$X$ by an open neighborhood of $x$, we may assume that $X$ is
defined in some $\AAA^N$ by a regular sequence. In this case
$I_{Z_1}={\rm Jac}_X$, and we deduce that 
$J_r=\OO_X$. Conversely, suppose that
$\OO_{X,x}$ is not locally complete intersection, 
and after restricting to an affine neighborhood of $x$, assume that
we have a closed embedding $X\subset A$ as in our general setting.
Note the by assumption, for every 
 complete intersection $M$ in $A$ that contains $X$, the ideal
$((I_M\colon I_X)+I_X)/I_X$ is contained in the ideal $\frmm_x$ defining $x\in X$. On the other hand, following the notation in the proof of
Corollary~\ref{cor3_appendix} 
we see that given
$L= (i_1,\ldots,i_c)$ with $1\leq i_1<\cdots<i_c\leq d$, and $I_M=(F_{i_1},\ldots, F_{i_c})$
we have  
$\frb_L\subseteq ((I_M\colon I_X)^r+I_X)/I_X\subseteq\frmm_x$. Therefore 
$\frb\subseteq\frmm_x$, and since $J_r$ and $\frb$ have the same support
(they even have the same integral closure), we conclude that $x$ lies in the support of $J_r$.
\end{remark}

\bibliographystyle{amsalpha}

\end{document}